\documentclass[preprint,12pt,journal=cmame]{elsarticle} 
\usepackage{amsfonts,amsmath,amssymb}			
\usepackage{blindtext}
\usepackage{bm}
\usepackage{color}
\usepackage{dcolumn}
\usepackage{fancyhdr}	
\usepackage{float}
\usepackage{graphicx}
\usepackage{hyperref}
\usepackage{inputenc}
\usepackage{natbib}
\usepackage{rotating}   
\usepackage{sidecap}
\usepackage{siunitx}											
\usepackage{textcomp}
\usepackage{textgreek}
\usepackage{upgreek}
\usepackage{amsmath}
\usepackage{amssymb}
\usepackage{graphbox,graphicx} 
\usepackage{verbatim}

\newcommand{\deal}{\texttt{deal.II} } 
\newcommand{\bemstokes}{\texttt{BEMStokes} }

\graphicspath{{figures/}}

\begin{document}

\title{MicroROM: An Efficient and Accurate Reduced Order Method to Solve Many-Query Problems in Micro-Motility}

\author[add1]{Nicola Giuliani}
\author[add1]{Martin W. Hess}
\author[add1,add2]{Antonio DeSimone}
\author[add1]{Gianluigi Rozza\thanks{grozza@sissa.it}}
\address[add1]{SISSA~---~International School for Advanced Studies,
 Via Bonomea 265, 34136 Trieste, Italy }
\address[add2]{Sant'Anna~---~Scuola Universitaria Superiore Pisa,
Piazza Martiri della Libert\`{a}, 56127 Pisa, Italy}
\ead{grozza@sissa.it}
    

\begin{abstract}

In the study of micro-swimmers, both artificial and biological ones, many-query problems arise naturally.
Even with the use of advanced high performance computing (HPC), it is not possible to solve this kind of problems in an acceptable amount of time. 
Various approximations of the Stokes equation have been considered in the past to ease such computational efforts but they 
introduce non-negligible errors that can easily make the solution of the problem inaccurate and unreliable.
Reduced order modeling solves this issue by taking advantage of a proper subdivision between a computationally expensive offline phase and a fast and efficient online stage.

This work presents the coupling of Boundary Element Method (BEM) and Reduced Basis (RB) Reduced
 Order Modeling (ROM) in two models of practical interest, obtaining accurate and reliable solutions to different many-query problems. 
Comparisons of standard reduced order modeling approaches in different simulation settings and a comparison to typical approximations to Stokes equations are also shown.
Different couplings between a solver based on a HPC boundary element method for micro-motility problems and reduced order models are presented in detail.
The methodology is tested on two different models: a robotic-bacterium-like  and an Eukaryotic-like swimmer, and in each case two resolution strategies for the swimming problem, the split and monolithic one, are used as starting points for the ROM. An efficient and accurate reconstruction of the performance of interest is achieved in both cases proving the effectiveness of our strategy.

\section*{Highlights}
\begin{itemize}
   \item We show how interesting micro-swimmer models can be cast as many-query problems, suitable for parameteric model reduction.
   \item A comparison to established simplified hydrodynamic models shows the superior accuracy of a ROM approach. 
   \item Complex many-query micro-swimming problems are solved efficiently by coupling reduced order modeling with an accurate boundary element method solver.
\end{itemize}

\end{abstract}

\begin{keyword}
Micro-motility, BEM, Reduced Order Modeling, Optimization, Many-query problems.
\end{keyword}

\maketitle

\section{Introduction and Motivation}
\label{sec:introduction}

Many interesting phenomena depend on the swimming behavior of motile cells, so there is a growing need of accurate, reliable and efficient computational methods to be applied in the study of 
micro-organisms swimming in a fluid. 
Understanding the behavior of micro-swimmers can give interesting insights on many complex biological processes: 
the spread of a pathogen~\cite{Josenhans2002}, the reproductive efficiency of sperm cells~\cite{Gray1955,Ishimoto2017a}, 
the ability to change a motility strategy depending on the environment~\cite{Noselli2019,Shum2015, Pimponi2016}, and many others.
An accurate and reliable simulation of the swimming mechanisms makes also possible the rational optimization and design of artificial micro-robots, 
mimicking such behaviors~\cite{Passov2012,AlougesDeSimoneZoppello,Dai2016,Keaveny2013,Walker2015,Gutman2016}.
Eukaryotic swimmers use complex flagellar beatings to achieve motion. 
For example \emph{Chlamydomonas Reinhardtii} exploits a couple of flagella that can beat symmetrically to produce an oscillating 
motion~\cite{Guasto2010,Klindt2015}, and \emph{Euglena Gracilis}  uses a single flagellum~\cite{Rossi2017,Tsang2018} executing a non planar beat.

Given the very small characteristic length scale of the problem, inertia is negligible and the fluid can be well approximated using the Stokes flow~\cite{Purcell1977,Purcell1997}, for which a vast variety of simulation tools exist, among the others we recall the shape optimization results of~\cite{Fumagalli2015,Ta2018}. Swimmers move both in bounded and unbounded domains and their swimming mechanisms usually involve very large geometrical deformations to achieve the net rigid motion~\cite{Purcell1977,Lauga2006}. This kinds of deformations make the use of standard Finite Element Method (FEM) very challenging. A possible alternative is the Immersed Boundary Method (IBM) which has been successfully applied to simulate micro-swimmers in bounded channels~\cite{Shi2016}. The simulation of swimmers in unbounded domains poses a challenge even for IBM.
The Boundary Element Method (BEM)~\cite{Steinbach2008} can deal with such large deformations efficiently if compared to other simulation techniques such as Finite Element or Finite Volume Methods.
While parallelized, high-performance solvers allow very accurate simulations, many-query problems like optimization and inverse problems in a repetitive computational environment still remain a challenge. 
Optimization procedures which address the different possible configurations and the reconstruction of the movements of a biological swimmer require to run 
many simulations quickly increasing the computational cost of the overall algorithm. 
Efficient BEM implementations based on open source High Performance Computing (HPC) libraries, are available~\cite{Giuliani2018,Giuliani2018SoRo,BEMStokesrepo}, 
but given the fact that the computational cost of BEMs scales at least quadratically with the problem complexity, the overall time required by many simulations quickly becomes unbearable even on modern computational architectures. A possible solution comes from rigorous mathematical approximations of the convolution integrals defining the BEM, among the others we recall the Fast Multipole Method~\cite{Greengard1987,Giuliani2018}, the Sparse Cardinal Sine Decomposition~\cite{Alouges2015,Alouges2017} and Reduced Order Modeling (ROM) which takes advantage of the classical many-queries settings employing a decomposition between a slow offline phase and a fast online evaluation, see~\cite{Rozza:ARCME}. We aim to solve  many-queries problems in micro-motility so we couple ROM with BEM.

In the present work we consider two different test-cases: the optimization of a bacterium-like artificial swimmer and the stroke reconstruction of an Eukaryotic-like organism. 
These examples present different challenges and they highlight the many-query scenarios that can be studied.
Many simplified models exist to ease the computational costs, such as Resistive Force Theory~\cite{Gray1955,Lighthill1976} or neglecting the hydrodynamic interactions of different parts of 
the swimmer~\cite{Purcell1997}, but the errors introduced are often not acceptable~\cite{Rodenborn2013, Giuliani2018SoRo}.

\begin{figure}
\centering
\includegraphics[width=.99\textwidth]{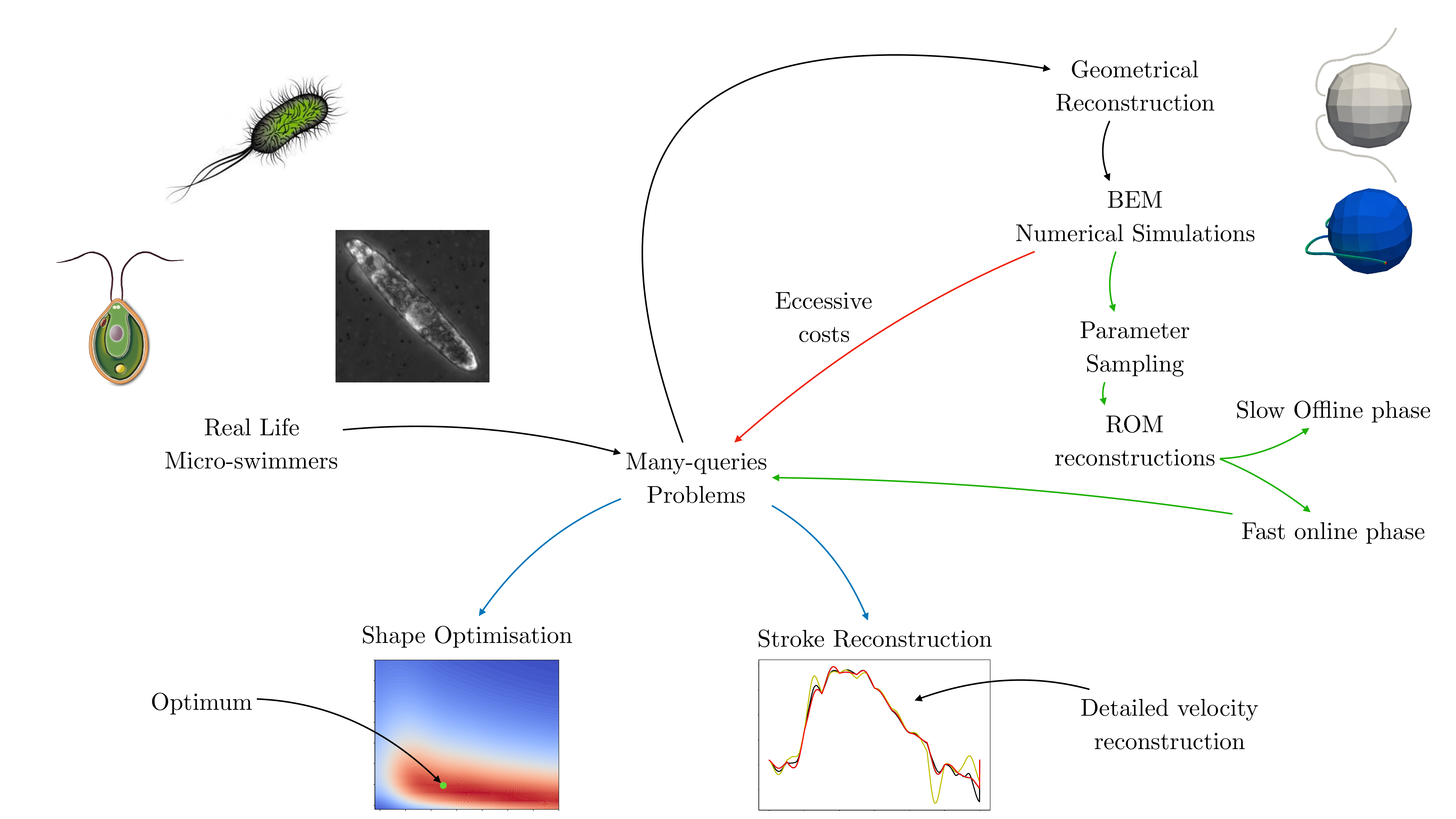}
\caption{The observation of real micro-swimmers as bacteria or more complex Eukaryotic species (as \emph{Chlamydomonas Reinhardtii} or \emph{Euglena Gracilis}) reveals the presence of different many-query problems. We work on two  of these problems: shape optimization and the detailed reconstruction of the swimming stroke. Even on modern computational architectures the direct numerical resolution of these  problems is unfeasible. We propose a different approach which couples Reduced Order Modeling and the Boundary Element Method. 
}
\label{sketchPaper}
\end{figure}

The motility of micro-swimmers can be understood using accurate simulations of the underlying fluid-dynamics. 
To enable many-query computations, reduced order methods are used upon a parallelized high-performance simulation software.
Many-query computations then allow to address simulation problems of larger scale, such as iterative optimization and inverse problems. Figure~\ref{sketchPaper} depicts the many query scenario coming from the observation of micro-swimmer and the solution we propose.

The coupling of BEMs and Reduced Order Modeling (ROM), see~\cite{Rozza:ARCME}, can reduce dramatically the computational costs of these many-query simulations. ROM have been recently coupled with Hybridisable discontinuous Galerkin Methods (HDG) to solve geometrically parametrized Stokes flows around micro-swimmers in~\cite{HuertaArxiv}. 
ROM and BEM have been coupled to study the airflow around objects in 2D~\cite{Manzoni2015} but, as far as the authors know, this is the first time these methods are combined to solve micro-swimmer 
problems on arbitrary geometries in 3D.
The model order reduction procedure relies on an offline-online decomposition, \cite{Rozza:ARCME}. 
During a computationally intensive offline phase, full-order solutions are sampled, an affine parameter-dependency of the 
system matrices is approximated and low-order quantities for the ROM are pre-computed. 
Due to the large computational load, the offline phase is typically done on a High Performance Computing  cluster. 
The online phase is then used for fast many-query evaluations with the reduced order model. 

In contrast to the standard reduced basis approach, we do not completely decouple all high-order dependencies from the reduced order model.
This is to allow a maximum of geometric flexibility of the micro-swimmer. In particular, the mesh is being generated at each new parameter value during the reduced order solve.
The alternative would be to approximate mesh movement in an affine fashion, but this would limit the possible application scenarios significantly. 
To reach the maximum possible flexibility the mesh is being generated with the software \emph{Blender}~\cite{BlenderOnlineCommunity2016}.
This provides only single precision accuracy, which limits a bit the maximum accuracy that the model reduction can attain. 
Nevertheless we believe the use of \emph{Blender} to be a strategic advantage since (i) it allows to generate complex arbitrary geometric parameterization in a straightforward way, (ii) 
it is an open source software, (iii) it has a large community of users.

The two standard approaches of reduced basis model reduction are explored, i.e., the proper orthogonal decomposition (POD) and the Greedy approach~\cite{Rozza:ARCME}.
The empirical interpolation method (EIM)~\cite{Barrault2004} is used to approximate an affine parameter dependency of the system matrices.
The focus is on the POD method, which is employed in both examples for all the full order resolution strategies. The Greedy approach is tested on one example to have a comparison with the POD available.

The Eukaryotic-like example exhibits a time-dependent boundary control.
The reduced order modeling of time-depending problems usually makes use of a combined POD-Greedy approach, i.e., 
POD in time and greedy in other parameters~\cite{Nguyen2009}, \cite{Hesthaven2016} and\cite{Haasdonk2008}. This is useful, as most time-dependent PDEs involve inertia, i.e., a time-derivative.
In the micro-swimming scenario, a time-derivative is not present and time can thus be treated as any other physical or geometrical 
parameter from the model reduction point of view.

As a all our approach produces very accurate results while achieving a drastic reduction in computational cost. We demonstrate this result on two different test cases of interest: the optimization of artificial micro-robots and the reconstruction of the stroke starting from a reduced set of geometrical configurations.

The organization of this work is as follows. Section 2 introduces the general boundary element framework in the continuous and discrete setting, while Section 3
presents two model swimmers mimicking real-life microorganisms. 
Section 4 explains the reduced modeling procedure and Section 5 presents detailed numerical results.
Section 6 summarizes, concludes and provides further perspectives.

\section{Problem Description}
\label{sec:model}
We aim to solve the swimming problem, i.e., to recover the swimming motion from a given history of shape changes. To do so, we define the mathematical and numerical modeling required to properly study micro-organisms swimming in a fluid medium. In particular, the swimming model is composed by a kinematic part describing the possible motions for the swimmer and a dynamic part defining the conservation laws for the fluid domain. The final resolution strategy couples these two models and is solved numerically either in a monolithic approach~\eqref{monolithic_system_complete_1}
 or in a split approach. 

\subsection{Swimming model}
Following \cite{Maso2011,DalMaso2015} the swimmer is represented as a time-dependent bounded open set $B_t \in \mathbb{R}^3$. The map $\chi:  \bar{B}_0 \subset \mathbb{R}^3 \times
[0,T] \rightarrow \mathbb{R}^3$ defines the position $x$ at time $t$ of a material point $X$ of the swimmer, namely
\begin{equation}
x(X,t) = \chi(X,t) = q(t) + R(t)s(X,t),
\label{GenericDispl}
\end{equation}
where $q(t)$, $R(t)$, $s(X,t)$ are the position of a point of the swimmer, the rotation of the system local frame and the shape at time $t$. The time derivative of \eqref{GenericDispl} defines the velocity of any material point on the swimmer,
\begin{equation}
\begin{aligned}
u_{s} = &\dot{x} = \dfrac{\partial \chi(X,t)}{\partial t} \\=&  \dfrac{d
  q}{d t} + R(t)\dfrac{\partial s(X,t)}{\partial t} + \dfrac{d
  R(t)}{d t} s(X,t)=\\
=&\dot{q}(t) + R(t)\dot{s}(X,t) + \omega(t) \wedge (R(t) s(X,t)).\\
\end{aligned}
\label{u_swimmer}
\end{equation}
We assume $s(X,t)$ to be known, we define ``shape velocity" the quantity $R(t)\dot{s}(X,t)$, therefore the unknowns are the rigid movements $q(t)$, $R(t)$, which define the linear and angular velocities $\dot{q}(t), \omega(t)$. 
In the case of self-propelled swimmers, only viscous drag is acting on the swimmer so the usual momentum balance laws read
\begin{subequations}
\begin{equation}
\int_{\Gamma} f(x) d\gamma(x)  =0,
\label{LinearMomentumBalanceDummy}
\end{equation}
\begin{equation}
\int_{\Gamma} f(x) \wedge (x-x_0) d\gamma(x) =0.   
\label{AngularMomentumBalanceDummy}
\end{equation}
\label{DummyMomentumSystem}
\end{subequations}
The viscous tractions $f$ are given by the action of the Cauchy stress tensor $\sigma$, see~\cite{Gurtin}, namely
\begin{equation}
f = \sigma(u,p) n,
\label{ForceContinuumMechanics}
\end{equation}
where $n$ is the outer unit normal vector to the surface, and $u,p$ represent the velocity and the pressure
in the fluid.

\subsection{Fluid model}
We study micro-swimmers immersed in a Newtonian incompressible fluid, 
by writing the non-dimensional formulation of the incompressible Navier-Stokes equation in a generic domain $\Omega$, see~\cite{Gurtin}, as
\begin{subequations}
\begin{equation}
\nabla \cdot u = 0 \quad  \text{ in } \Omega  \
\label{MassConservationNStokesintro}
\end{equation}
\begin{equation}
Re \left( \bar{\sigma}\frac{\partial u}{\partial t} + (\nabla u) u \right)= \Delta u -\nabla p \quad \text{ in } \Omega,
\label{LinearMomentumBalanceNStokesintro}
\end{equation}
\label{NavierStokesSystemintro}
\end{subequations}
where $u, p$ are non dimensional velocity and pressure and $Re$ is the Reynolds number defined as
\begin{equation}
Re = \frac{\rho U L}{\mu},
\label{Reynolds}
\end{equation}
while
\begin{equation}
\bar{\sigma} = \frac{\omega L}{U}
\label{Womersley}
\end{equation}
is the Womersley number.
We consider the typical length, velocity and frequency for a micro-swimmer as $L = 10^{-5} m$, $U = 10^{-5} m/sec$ and $\omega = 10^{2} Hz$, we obtain $Re=10^{-4}$ and $\bar{\sigma}Re = 10^{-2}$. Therefore the system of equations governing the flow is the Stokes system,
\begin{subequations}
\begin{equation}
\nabla \cdot u = 0 \quad \text{ in } \Omega  \
\label{MassConservationStokesintro}
\end{equation}
\begin{equation}
\Delta u -\nabla p = 0 \quad \text{ in } \Omega.
\label{LinearMomentumBalanceStokesintro}
\end{equation}
\label{StokesSystemintro}
\end{subequations}
We identify the boundary of the fluid domain with the boundary of the swimmer $\partial \Omega = \Gamma$, we write Dirichlet's boundary conditions for \eqref{StokesSystemintro} as
\begin{equation}
u = u_{s},
\end{equation}
where $u_s$ denotes the velocity of the swimmer on the boundary defined in \eqref{u_swimmer}. 
The system~\eqref{StokesSystemintro} with Dirichlet boundary conditions is well posed and admits a unique solution $(u,p) \in ((H^1(\Omega))^3, L_2(\Omega))$.
Following \cite{Pozrikidis1992} we rewrite \eqref{StokesSystemintro} using the fundamental solution to get the representation formula for the Stokes system
\begin{equation}
u_i(x) - \int_{\Gamma} W_{ijk}(x,y) n_k(y) u_j(y) d\gamma_y = \int_{\Gamma}
G_{ij}(x,y) f_j(y)d\gamma_y \ \forall x \in \mathbb{R}^d
  \setminus \Gamma,
\label{StokesRepGeneral}
\end{equation}
where $G, W$ are the first two fundamental solution for the Stokes system.
We take the trace of \eqref{StokesRepGeneral} to obtain the Boundary Integral Equation (BIE) of the Stokes system as
\begin{equation}
\alpha(x)u_i(x) - \int_{\Gamma}^{PV} T_{ijk}(x,y) n_k(y) u_j(y) d\gamma_y = \int_{\Gamma}
G_{ij}(x,y) f_j(y)d\gamma_y \ \forall x \in \Gamma,
\label{StokesBIEGeneral}
\end{equation}
where the integral on the left is computed in the Cauchy principal value
sense, and $\alpha$ represents its Cauchy principal value. 
We define the single and double layer operator $V, H$ as
\begin{subequations}
\begin{equation}
H: (H^{\frac{1}{2}}(\Gamma))^3 \rightarrow (H^{\frac{1}{2}}(\Gamma))^3,
\label{doublelayerfunctional}
\end{equation}
\begin{equation}
V: (H^{-\frac{1}{2}}(\Gamma))^3 \rightarrow (H^{\frac{1}{2}}(\Gamma))^3,
\label{singlelayerfunctional}
\end{equation}
\label{boundaryoperatorsfunctional}
\end{subequations}
where $ (H^{\frac{1}{2}}(\Gamma))^3$ represents the space of the traces of the functional $(H^1(\Omega))^3$ and $ (H^{-\frac{1}{2}}(\Gamma))^3$ is its dual space.
We rewrite \eqref{StokesBIEGeneral} using the operators~\eqref{doublelayerfunctional} and~\eqref{singlelayerfunctional} as
\begin{equation}
\left[ \alpha I - H \right] u = -K  u=-  V  f.
\label{StokesOperator}
\end{equation}
We remark that we use the complete traction $f$ as given by~\eqref{ForceContinuumMechanics}. For this reason, differently from standard Finite Element simulations of~\eqref{NavierStokesSystemintro}, we don't need to solve for the pressure $p$ and the problem is not a saddle point problem. The usage of the fundamental solutions $G, T$ guarantees that the velocity in the fluid domain is automatically divergence free, see~\cite{Steinbach2008} for greater details. For these reasons operator $V$ is coercive, similarly to what happens to the classic Stokes system~\eqref{NavierStokesSystemintro}  if we consider a velocity domain which is divergence free.

\subsection{Numerical resolution of the swimming problem}
\label{sec:fom_resolution}
To solve the swimming problem means to recover the translation $q(t)$ and rotation $R(t)$ given the history of shape changes. 
For this purpose we use the numerical methodology presented and validated in~\cite{Giuliani2018SoRo} to discretize \eqref{StokesBIEGeneral} into a Boundary Element Method (BEM) \footnote{In particular we use the open source software \bemstokes which is freely available under LGPL license v.2.1 on github~\cite{BEMStokesrepo}}. We discretize both the geometry of the swimmer and the unknowns $f, u$ using  Lagrangian finite element and we apply a collocation method to derive the BEM linear system
\begin{equation}
\hat{K}  \hat{u}= \hat{V} \hat{f},
\label{StokesOperator2}
\end{equation}
where $\hat{\dots}$ represents the finite dimensional counterpart of the operators and functionals represented in~\eqref{StokesOperator}. Following~\cite{Steinbach2008} we express the tractions as a function of the velocity introducing the discretized Dirichlet to Neumann map 
\begin{equation}
\hat{T} = \left[ \hat{V}^{-1} \hat{K} \right],
\label{DNMap}
\end{equation}
so that $\hat{f} = \hat{T}\hat{u}$. We will drop the $\hat{\dots}$ in the rest of the paper for the sake of simplicity.
Following~\cite{Giuliani2018SoRo} we rewrite \eqref{u_swimmer} using a suitable set of basis functions for the rigid velocity $P(X,t)$ and introducing the shape velocity in the swimmer body frame $v(X,t) = R(t) \dot{s}(X,t)$, as
\begin{equation}
\begin{aligned}
u_{s} &= \dot{q}(t) + \omega(t) \wedge R(t) s(X,t) + v(X,t) \\
&=\sum_{i=1}^{N_{r}} p_i(X,t) \dot{p}_i(t)+ v(X,t) \\
&=P(X,t) \dot{p}(t) + v(X,t).
\end{aligned}
\label{u_swimmer2}
\end{equation}
We assume $\dot s(X,t)$ to be known, so at time $t$, given $R(t)$, we compute $v(X,t)$ and we write 
\begin{subequations}
\begin{equation}
f_{shape} = T v,
\label{shapeforce}
\end{equation}
together with
\begin{equation}
F_{rigid} = T P.
\label{rigidforce}
\end{equation}
\label{DNforces}
\end{subequations}
We obtain the rigid velocity coefficients $\dot{p} (t)$ using~\eqref{DummyMomentumSystem} as
\begin{equation}
P^T M F_{rigid} \dot{p} + P^T M f_{shape} = 0,
\label{numericalmomentumbalance}
\end{equation}
where $M$ is the mass matrix taking care of the surface integrals of~\eqref{DummyMomentumSystem}.
We then solve ~\eqref{numericalmomentumbalance} to get 
\begin{equation}
\dot{p} = - \left[ P^T M F_{rigid}\right]^{-1} P^T M f_{shape}.
\label{solvingequation}
\end{equation}
We remark that $ \left[ P^T M F_{rigid}\right]$ is a symmetric negative definite matrix called grand resistive matrix, so its inversion in~\eqref{solvingequation} is legit. We call this resolution strategy, where we compute the tractions associated to any possible velocity and then combine them, the \emph{split approach}.

Another possible resolution strategy is the so-called \emph{monolithic approach}~\cite{Pimponi2016,Giuliani2018SoRo}, which is based on a single resolution to obtain the overall traction $f$ together with the rigid velocity coefficients. We write the discretized BIE together with the constraint~\eqref{numericalmomentumbalance} as
\begin{equation}
\begin{bmatrix}
[V] & -[K]P \\
P^T [M]  & 0
\end{bmatrix}
\begin{bmatrix}
f\\
 \dot{p}
\end{bmatrix}
=
 \begin{bmatrix}
[K] v\\
 0
\end{bmatrix},
\label{monolithic_system_complete_1}
\end{equation}
the solution of~\eqref{monolithic_system_complete_1} is equivalent to the split approach. The differences lie in the fact that the monolithic approach only needs one resolution while the split needs seven. Moreover even the spectral properties of the matrices change if we impose the constraint during the BEM resolution. We now apply ROM to both full order approaches to highlight pros and cons of the two resolution strategies.

Standard numerical integration of $\dot{p}(t)$ provides the rigid displacements and rotations. 
The \bemstokes library efficiently combines several algorithms and ideas in a flexible, modular, and extendible way. 
It exploits distributed memory parallelism (MPI) using an automatic splitting of the workload at algebraic level using on the graph partitioning tool METIS~\cite{GeorgeKarypis1998}, 
with the high performance computing libraries Trilinos~\cite{Heroux2005}, and \deal~\cite{dealII90}, used to  tackle distributed linear algebra. 
A similar combination has been successfully applied to achieve high computational efficiency in fluid dynamics, as demonstrated in ASPECT~\cite{Kronbichler2012} and $\pi$-\texttt{BEM}~\cite{Giuliani2018}.

\section{Application: micro-swimmers test cases}
\label{sec:micro_robot_models}

We consider two models inspired directly from microswimmers found in nature. 
At the micron scale, there exists a wide variety of organisms \cite{harris1989the,archibald2017handbook,Guasto2010}, which swim using the motion of flagella. 
The resulting flow can be well modelled using Stokes inertia-free equations~\cite{Purcell1977,Lauga2006,Ishimoto2017}. 
The mathematical model swimmers presented in this section are inspired one from a bacterium and one from the Eukaryotic swimmer \emph{Chlamydomonas Reinhardtii}. 
The key difference lies in the motion of the flagellum: while bacteria, sketched in Figure~\ref{BacteriumSetting} can only rotate the root of the flagellum, which is otherwise a passive structure, see~\cite{Son2013,Shum2015}, 
Eukaryotic swimmers can exploit a more complex flagellar architecture, called axoneme, to control the shape of the flagellum~\cite{Porter2000}. We sketch these different configuration in Figure~\ref{chlamy}.  

In particular, for the bacterium presented in Section~\ref{subsec:micro_robot_models_bacterium} we consider two geometrical parameters, the head radius and the number of flagella windings.
The Eukaryotic swimmer of Section~\ref{subsec:micro_robot_models_chlamydomonas} undergoes a complex flagella movement in a parametrized time trajectory, i.e., a time-dependent 
boundary control is considered.

\subsection{Bacterium-like microswimmer}
\label{subsec:micro_robot_models_bacterium} 

\begin{figure}
\centering
\includegraphics[width=.48\textwidth]{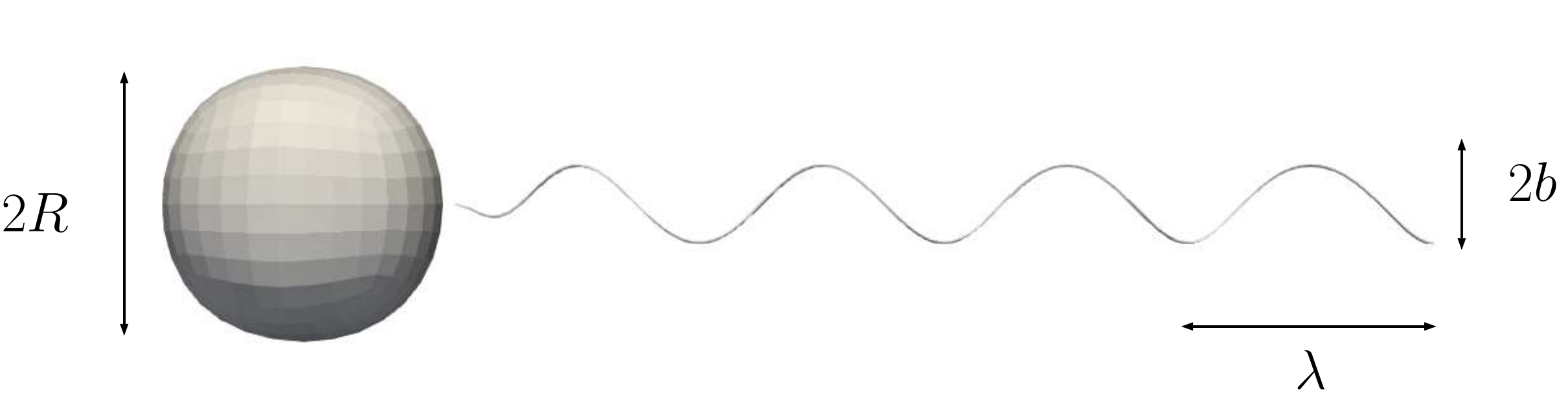}
\caption{Sketch representing a simple robotic-like bacterium composed by a spherical head of radius $R$ and a helical tail of pitch $\lambda$ and amplitude $b$.}
\label{BacteriumSetting}
\end{figure}

Following~\cite{FUJITA2001, Phan-Thien1987, Shum2010, Pimponi2016} we prescribe a robotic micro-swimmer to be composed by a rigid spherical head of radius $R$ 
and a helical tail of width $b$ which is rotating with respect to the head with constant velocity $\omega$. 
The centerline of the tail is given by 
\begin{equation}
{r} = (x,y,z) = (x, b E(x) cos(kx-\omega t), b E(x)
sin(kx-\omega t)),
\label{spiralThien}
\end{equation}
with
\begin{equation}
E(x) = 1 - e^{-(k_E x)^2}.
\end{equation} 
Figure~\ref{BacteriumSetting} depicts the micro-swimmer under consideration.
Following~\cite{FUJITA2001} we assume the pitch of the helix 
$\lambda = 2 \pi / k$, $k_E = k$, the flagellum has a total number of turns
of $N_\lambda$ and the width of the flagellum is $b=\lambda/(2 \pi)$. The flagellum thickness is given by 
\begin{equation}
d =  0.02  \frac{\pi}{4}  R.
\label{thickness}
\end{equation}
We discretize both the geometry and the unknown functional of~\eqref{StokesOperator2} using linear Lagrangian finite elements for an overall number of $N_\delta=2430$ degrees of freedom. 
We remark that we are considering organisms swimming in free space, in this case a rotational symmetry exists during the stroke (the rotation of the tail w.r.t the head), 
therefore we simulate the complete stroke of the bacterium-like model solving a single time instant and then we rotate the corresponding instantaneous results to get the results of a complete stroke.

We consider a two-dimensional parameter domain, with varying number of helix turns $N_{\lambda} \in [0.4, 4.0]$ and head radius $R_{head} \in [0.4, 4.0]$, so the 
parameter domain is $\mathcal{P} = [0.4, 4.0]^2$.
Following \cite{FUJITA2001} we use~\eqref{thickness} to let the thickness of the flagellum increase proportionally together with the head radius.  
We aim to use ROM to find the optimal parameter set $\boldsymbol{\mu} \in \mathcal{P}$ maximizing a chosen performance measure. Following~\cite{Purcell1997,Phan-Thien1987} we consider an 
energetic efficiency, called Lighthill efficiency. Following~\cite{FUJITA2001} we compute the axial velocity $U_{axial}$ as the projection of the velocity along the angular velocity of the tail, we write
\begin{equation}
U_{axial} = {\dot{q}} \cdot \frac{{\omega} - {\Omega}}{|{\omega} - {\Omega}|},
\label{U_axial}
\end{equation}
where $\dot{q}, {\Omega}$ are the rigid velocities of the swimmer and ${\omega}$ is the relative angular velocity between head and tail.
Using~\eqref{U_axial} we write the Lighthill efficiency as
\begin{equation}
\eta_{Lighthill}=\frac{\text{effective power}}{\text{input power}} =\frac{D_{head} U_{axial}}{T_{motor} \omega} = \frac{K_{head} U_{axial}^2}{T_{motor} \omega},
\label{eta_lighthill}
\end{equation}
where $K_{head}$ is the drag coefficient for the head considered alone in free space.
We highlight that~\eqref{eta_lighthill} is the ratio between the power expanded to move the payload (the head) and the overall power expanded by the swimmer.

\subsection{Eukaryotic-like microswimmer}
\label{subsec:micro_robot_models_chlamydomonas}

\begin{figure}
\centering
\includegraphics[width=.3\textwidth]{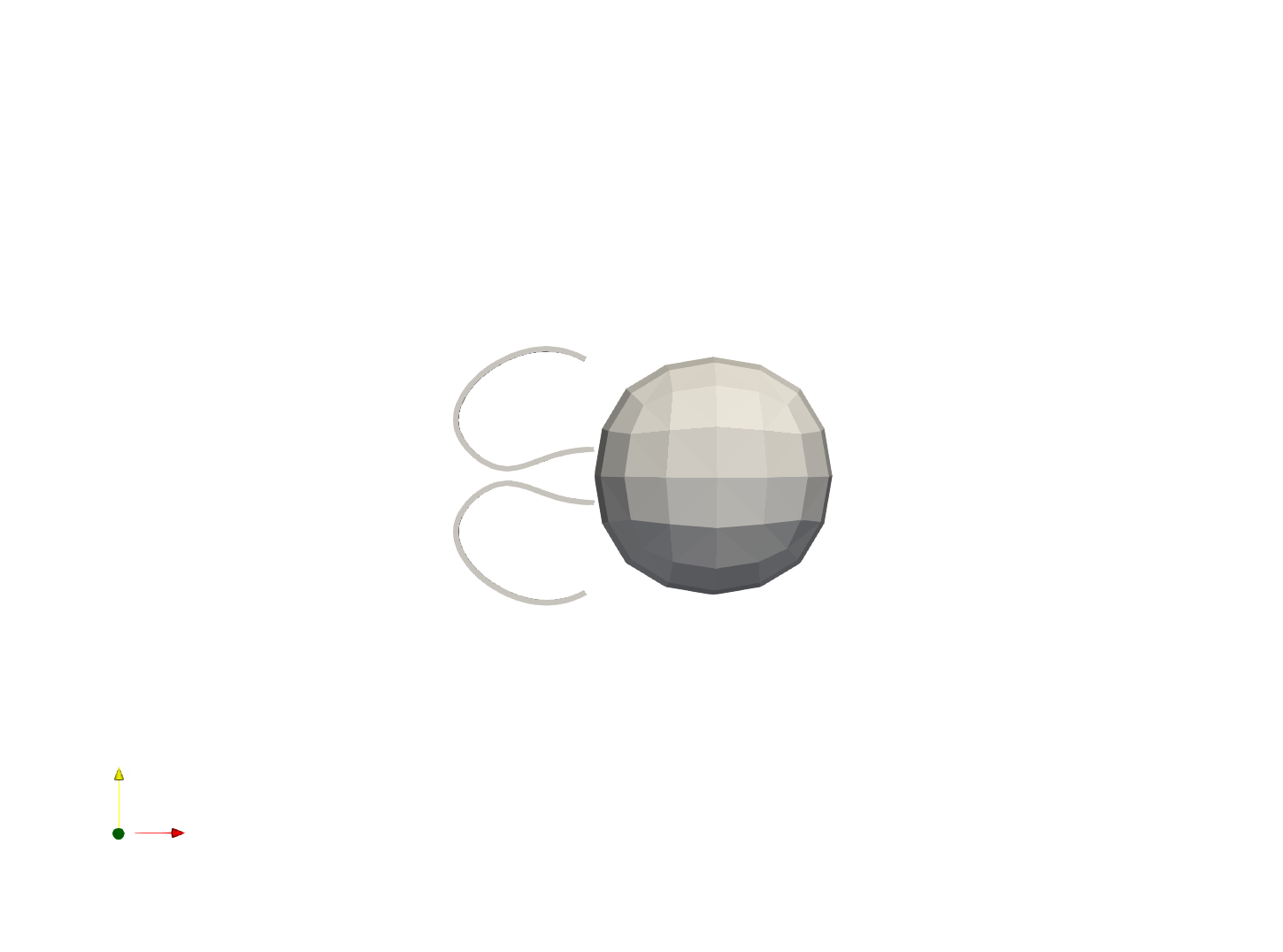}\includegraphics[width=.3\textwidth]{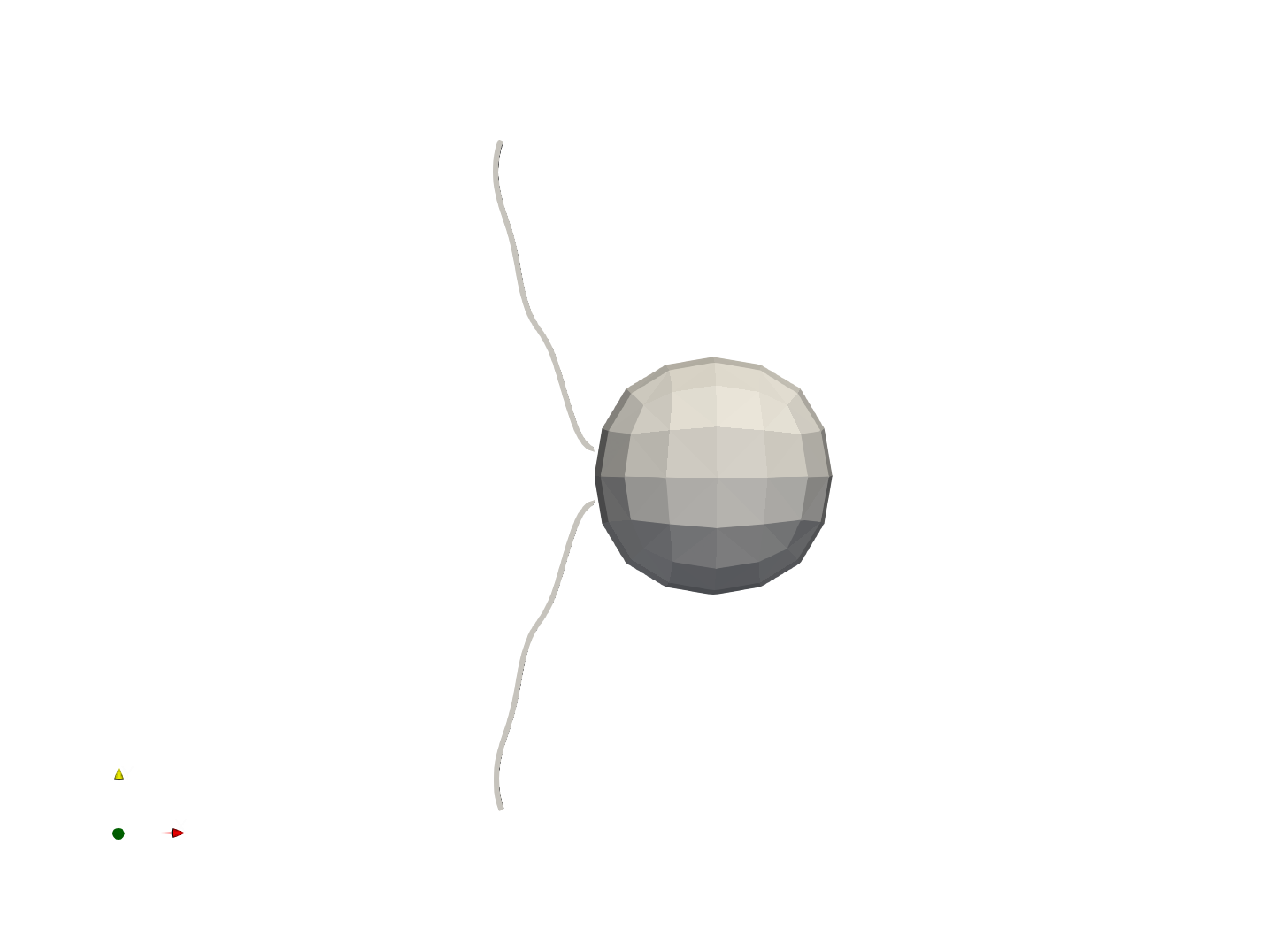}\includegraphics[width=.3\textwidth]{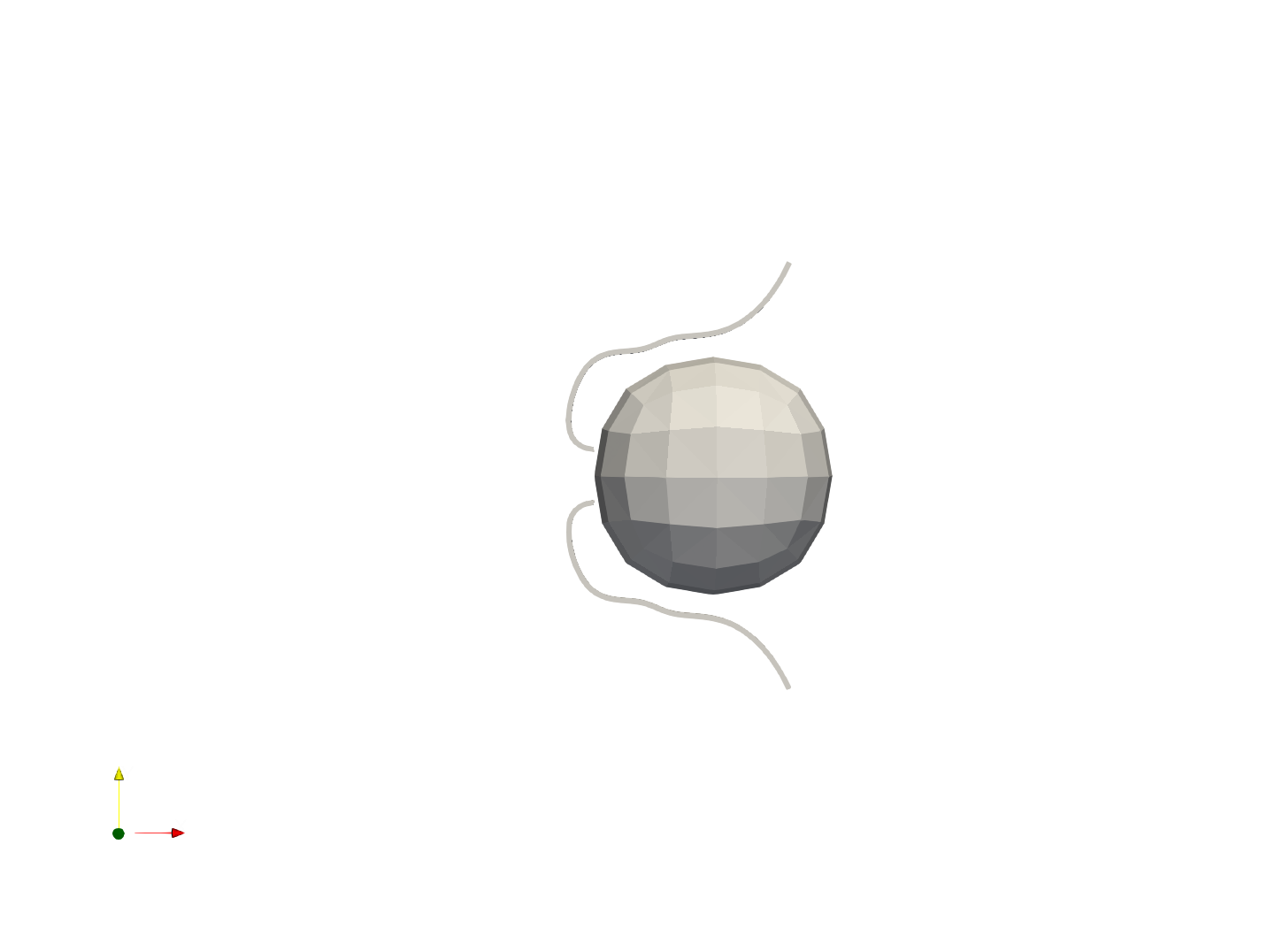}
\caption{Sketch of an Eukaryotic swimmer representing \emph{Chlamydomonas Reinhardtii}. 
We consider 1200 frames to represent the complete stroke: on the left frame 0, on the center frame 400 and on the right frame 800. 
}
\label{chlamy}
\end{figure}

To prove the versatility of our proposed methodology we apply reduced order modeling to a different kind of swimmer, namely we consider an Eukaryotic swimmer mimicking a \emph{Chlamydomonas Reinhardtii} specimen. 
This organism has been extensively studied in recent years~\cite{Drescher2010a,Guasto2010,Geyer2013,Drescher2009} and it is composed by two flagella and a spherical cellular body. 
The flagella are activated in different ways thanks to the axoneme structure~\cite{Lauga2009} and we focus on the symmetric beating of the two cilia producing a ``breast stroke" like movement in 
the swimmer.  We have interpolated the experimental observation of~\cite{Guasto2010} using a standard cubic spline interpolation to obtain a continuous synthetic stroke. For further details on the numerical procedure the reader is referred to~\cite{GiulianiPHD}.
We parametrize the flagellum centerlines as two NURBS beating symmetrically with respect to the body centerline. We represent the periodic stroke using 1200 equi-distanced frames $\phi$, we plot the 
numerical mesh at three different times during the stroke in Figure~\ref{chlamy}. We consider a spherical cell body having diameter of $5 \mu m$ with two symmetrical beating flagella with length $8 \mu m$ 
and thickness $0.1 \mu m$, we use standard linear Lagrangian finite elements for an overall total of $N_\delta = 1902$ degrees of freedom.
This problem is a time-dependent boundary control problem and we remark that since the equations of motion~\eqref{StokesBIEGeneral} do not have any time dependence we can consider the frame number as a 
geometric parameter $\boldsymbol{\mu}$, indexing the current shape. 
We aim to reconstruct the complete history of rigid body velocities, and consequently the fluid flow around the swimmer, depending on the shape changes, in an efficient and reliable way, exploiting a ROM built from the knowledge of a set of frames. This is also a typical problem from an experimental point of view since the usual frame rates only allows for the reconstruction a discrete subset of shape changes during the stroke.

\section{Reduced Order Models for microswimmers}
\label{sec:rom}

A reduced order model (ROM) is a low-dimensional surrogate model, which approximates the original high-dimensional model over a parameter range of interest.
In the following, parameter-dependent quantities are indicated with the parameter vector $\boldsymbol{\mu}$.
The ROMs are created by means of proper orthogonal decomposition (POD) and reduced basis (RB) methods, see \cite{Rozza:ARCME} for a general overview and \cite{Lassila2014} for an overview with 
a focus on fluid dynamics. 
The computational speed-up is achieved through an offline-online decomposition, i.e., a compute-intensive offline phase generates a ROM, while a fast online phase is used for many-query 
ROM evaluations for optimization purposes in a repetitive computational environment. A pre-requisite for the offline-online decomposition is an affine parameter dependency. 
Since the parameter-dependencies considered here are not affine, an affine parameter dependency is approximated by means of \emph{empirical interpolation method} (EIM) \cite{Barrault2004}, \cite{Chaturantabut2010}. 
A projection space is determined to project the high-dimensional equations onto the ROM. 
During the online phase a low-order approximate solution is computed, independent of the high-order discretization size.
Two algorithms are used to compute the projection space: the POD and a residual-based greedy sampling.
Both approaches are established in the area of ROM for fluid-structure interactions, see \cite{10.3389/fams.2018.00018}, \cite{Lassila_2012}.

\subsection{Offline-Online decomposition}
\label{sec:offline-online}
Conceptually, reduced order modeling uses an offline-online decomposition, in which a compute-intensive offline phase determines an accurate
reduced order model (ROM), while the online phase performs fast evaluations. This is beneficial in a many-query or real-time context, since 
the ROM can be evaluated at low computational cost and accurately recovers the high-dimensional solution.

After applying the \emph{empirical interpolation method} in the geometry parameter, the parameter dependency is cast in an affine form.
Therefore, there exists an affine expansion of the matrices $V(\boldsymbol{\mu})$ and $K(\boldsymbol{\mu})$ in the parameter $\boldsymbol{\mu}$.
To achieve fast reduced order solves, the offline-online decomposition computes the parameter-independent projections offline, which are
stored as small-sized matrices of the order $N \times N$, with $N$ depending on the velocity and traction projection spaces. 
When solving for a new parameter online, the affine form \eqref{V_K_affine_expanded} is evaluated and the reduced order solution \eqref{ROM:final_proj} computed.

In contrast to the typical ROM approach, the dependence on the high-order discretization size is not fully removed in the online phase.
The mesh points are generated for the full order problem at each new parameter but only the mesh points relevant for the EIM are taken into account for further computations.
This allows to have a maximum of geometric flexibility, needed for these problems, as affinely parametrizing the mesh points would be only possible for special cases of movement. Moreover we do not apply any reduction in the resolution of the momentum balances~\eqref{numericalmomentumbalance} since the computation of the full order method is already efficient, therefore we focus the reduction strategy on the discretization of the BIE~\eqref{StokesOperator2}.

\subsection{Empirical interpolation method}
\label{sec:DEIM}

Since parametric variations in geometry affect the BEM system \eqref{StokesOperator} in a nonlinear way,  \emph{empirical interpolation method} (EIM) (\cite{Barrault2004}, \cite{Chaturantabut2010}) 
computes an approximate affine parameter dependency.
The matrix discrete empirical interpolation \cite{Negri2015} computes the decompositions

\begin{equation}
\begin{aligned}
 \sum_{i=1}^{Q_V} \Theta_V^i(\boldsymbol\mu) V_i  \approx V(\boldsymbol\mu), 
\quad \quad
  \sum_{i=1}^{Q_K} \Theta_K^i(\boldsymbol\mu) K_i \approx K(\boldsymbol\mu).
\label{V_K_affine_expanded}
\end{aligned}
\end{equation}

\noindent with scalar parameter-dependent coefficient functions $\Theta_V^i(\boldsymbol\mu)$ and $\Theta_K^i(\boldsymbol\mu)$ and parameter-independent matrices $V_i$ and $K_i$.
In particular, each coefficient function corresponds to a single matrix entry of $V$ or $K$, respectively.
Since the assembly of only a few (here $Q_a < 260$) matrix entries can be implemented efficiently, an approximation of $V$ and $K$ is readily available for each new $\boldsymbol\mu$.
As in the determination of the projection spaces $U_u, U_f$, at least $99.99\%$ of the POD energy is used to approximate the system matrices from the collected matrices during snapshot computation.
The EIM is performed for each of the two BEM matrices $V,K$ and the reduced mode are stored considering the $6$ sets of rigid velocities and the 
shape velocity for an overall total of $14$ reduced matrices. The affine expansions on the other hand, are the same for each rigid mode and the shape velocity.

\subsection{Proper Orthogonal Decomposition}

The POD samples uniformly distributed Stokes solutions $u, f$ over the parameter domain. The solutions are often called \emph{snapshots} in reduced order modeling.
A \emph{singular value decomposition} of the snapshots is computed and the most dominant modes are chosen as the projection space.
The most dominant modes corresponding to/ at least, $99.99\%$ of the POD energy form the projection matrix $U \in \mathbb{R}^{N_\delta \times N}$ and implicitly define the low-order space $V_N = \text{span} (U)$ 
and set the reduced model size $N$.

\subsection{Greedy Sampling}
\label{sec:greedy_split}
The greedy sampling builds the projection space iteratively from snapshot solutions. 
The parameter location of the next snapshot is chosen as the maximum of a residual-based error indicator, see  \cite{Rozza:ARCME} for details.
In particular we implement a classic residual-based error indicator to choose the snapshots for the force samples selection, while we use a projection-based error indicator for the selection of the velocity samples. The choice of different samples for all the possible snapshots allows for a greater accuracy and faster ROM computation.
The decay of the maximum residual over the iterations determines a stopping criterion of the greedy sampling.
As in the POD case, a projection matrix $U \in \mathbb{R}^{N_\delta \times N}$ is obtained after orthonormalizing the snapshots.

\subsection{Projection for the split approach}
\label{sec:split_rom}
In fact, two projection spaces are built, namely $U_u$ from sampling ${u}$ in \eqref{StokesOperator2} and $U_f$ from sampling $f$ in \eqref{shapeforce}.
The large-scale system \eqref{StokesOperator2} is then projected onto the reduced order space:
\begin{eqnarray}
U_f^T [{K} ](\boldsymbol{\mu}) U_u  U_u^T {u}(\boldsymbol{\mu}) = U_f^T [ {V} ](\boldsymbol{\mu}) U_f U_f^T {f}(\boldsymbol{\mu}) .
\label{ROM:final_proj}
\end{eqnarray}
The low order solutions ${u}_N(\boldsymbol{\mu}), {f}_N(\boldsymbol{\mu})$ approximates the large-scale solution as ${u}(\boldsymbol{\mu}) \approx U_u {u}_N(\boldsymbol{\mu}), {f}(\boldsymbol{\mu}) \approx U_f {f}_N(\boldsymbol{\mu})$.
In particular, the dimensions of the projection spaces $U_u$ and $U_f$ can be chosen independently from each other for every momentum and rigid velocity leading to the creation of $7$ different reduced order systems as the one depicted in~\eqref{ROM:final_proj}.
Choosing them independently has two major advantages: 
in case of the greedy sampling, less snapshot solutions need to be computed, and the evaluation of the ROM becomes faster since the ROM dimensions are smaller.
Once we solve the $7$ different ROM systems we can use~\eqref{DummyMomentumSystem} to obtain the six rigid velocities, in particular we use the full order mass matrix $M(\boldsymbol{\mu})$ together with the rigid modes to obtain
\begin{eqnarray}
 P^T(\boldsymbol{\mu}) [M(\boldsymbol{\mu})] \left[ U_f^1 f_N^1(\boldsymbol{\mu})|\dots|U_f^6 f_N^6(\boldsymbol{\mu})\right] p(\boldsymbol{\mu}) = -P^T(\boldsymbol{\mu}) [M(\boldsymbol{\mu})] U_f^{s} {f}_N^{s}(\boldsymbol{\mu})
\label{ROM:final_balances}
\end{eqnarray}
where $P(\boldsymbol{\mu}), M(\boldsymbol{\mu})$ represent the full order rigid modes and mass matrix respectively. The assembling time of the full order mass matrix is very fast given its sparsity with respect to the one required by the two BEM operators, so we choose to use the full order operator to have the maximum possible accuracy.


\subsection{Projection for the monolithic approach}
\label{sec:monolithic_rom}
The full order model presented in Section~\ref{sec:fom_resolution} is based on seven different application of the Dirichlet to Neumann map $T$ to retrieve the stress fields associated to shape and rigid velocities and then combine them to fulfill ~\eqref{DummyMomentumSystem} and retrieve the actual rigid velocity coefficients. Alternatively, we use the ``so-called'' monolithic approach~\cite{Pimponi2016,Giuliani2018SoRo} that solves the Stokes system imposing the constraints specified by the momentum balances~\eqref{DummyMomentumSystem} at the same time. We write the monolithic system as  
\begin{equation}
\begin{bmatrix}
[V](\boldsymbol{\mu}) & -[K](\boldsymbol{\mu})P(\boldsymbol{\mu}) \\
P^T(\boldsymbol{\mu}) [M](\boldsymbol{\mu})  & 0
\end{bmatrix}
\begin{bmatrix}
f(\boldsymbol{\mu})\\
 \dot{p}(\boldsymbol{\mu})
\end{bmatrix}
=
 \begin{bmatrix}
[K](\boldsymbol{\mu}) v(\boldsymbol{\mu})\\
 0
\end{bmatrix},
\label{monolithic_system_complete}
\end{equation}
where $[V](\boldsymbol{\mu}), [K](\boldsymbol{\mu}) $ are the Stokes operators defined in Section~\ref{sec:fom_resolution} and $[M](\boldsymbol{\mu})$ is the Mass matrix that takes care of the surface integration described in~\eqref{DummyMomentumSystem}, $f(\boldsymbol{\mu})$ represents the tractions on the boundary of the swimmer, $\dot{p}(\boldsymbol{\mu})$ are the rigid velocity coefficients while $v(\boldsymbol{\mu})$ is the shape velocity. We remark that the system~\eqref{monolithic_system_complete} consists of $N_\delta + 6$ equation in $N_\delta + 6$ unknowns. We repeat the procedure described in Section\ref{sec:offline-online} and we introduce the projection spaces $U_f, U_u$ obtaining
\begin{equation}
\begin{bmatrix}
U_f^T [V](\boldsymbol{\mu}) U_f & -U_f^T [K](\boldsymbol{\mu}) U_u {P}_N(\boldsymbol{\mu}) \\
P^T(\boldsymbol{\mu}) [M](\boldsymbol{\mu}) U_f & 0
\end{bmatrix}
\begin{bmatrix}
\tilde{f}(\boldsymbol{\mu})\\
 \dot{p}(\boldsymbol{\mu})
\end{bmatrix}
=
 \begin{bmatrix}
U_f^T [K](\boldsymbol{\mu}) U_u  {v}_N((\boldsymbol{\mu}))\\
 0
\end{bmatrix},
\label{monolithic_system_ROM}
\end{equation}
where ${f}_N(\boldsymbol{\mu}) = U_f^T f(\boldsymbol{\mu})$ is the reduced order representation of the tractions, ${P}_N(\boldsymbol{\mu}) = U_u^T P(\boldsymbol{\mu})$ is the projection of the rigid velocities and ${v}_N(\boldsymbol{\mu}) = U_u^T v(\boldsymbol{\mu})$ is the projection of the shape velocity. The system~\eqref{monolithic_system_ROM} consists of $N+6$ equation in $N+6$ unknowns, we remark that $[{V}_N](\boldsymbol{\mu})= U_f^T [V](\boldsymbol{\mu}) U_f$ and $[{K}_N](\boldsymbol{\mu}) = U_f^T [K](\boldsymbol{\mu}) U_u$ are computed using the empirical interpolation method described in Section~\ref{sec:DEIM}. 
The spaces $U_f$ and $U_u$ are obtained using a POD of samplings for traction and velocity field, the choice of the samplings has a key-role in determining the stability of the reduced monolithic system. We do not address this issue in detail in the present work, but, judging from our numerical experiments, we believe the best choice to be samplings of all the possible tractions and velocities the system experiences, namely $U_f$ is obtained from the POD of the shape traction $f_{shape}(\boldsymbol{\mu})$ and the rigid tractions $F_{rigid}(\boldsymbol{\mu})$ while $U_u$ comes from the sampling of both $v(\boldsymbol{\mu})$ and $P(\boldsymbol{\mu})$.

\subsection{Expected ROM accuracy}

Since the modeling and meshing is performed with the software \emph{Blender} v2.79~\cite{BlenderOnlineCommunity2016}, the 
input data for numerical simulation and ROM simulations are limited to single precision. 
We believe the choice of a wide-spread well documented open-source software to be of capital importance in order to guarantee flexibility to the micro-motility solver.
It is expected that there is in each step of (i) empirical interpolation, (ii) projection and (iii) solving in ROM space, a slight degradation in the accuracy.
We thus expect an accuracy between ROM and FOM of about four digits and an accuracy in an output quantity of about three to four digits.
Typically, the input data are double precision and a higher accuracy between ROM and FOM can be expected.

\subsection{Optimization using ROM, the coarse-fine approach}
The task is to optimize the swimming performance $\eta_{Lighthill}$ using a reduced order model built considering the head radius $(R_{head}$ and number of windings $N_{\lambda})$ over a two-dimensional parameter domain $\mathcal{P}$. 
To obtain an optimal speed-up, we propose a two-step method.
In the first step, a coarse sampling of the parameter domain is used to obtain the approximate location of the maximum swimming efficiency $\eta_{Lighthill}$. 
In the second step, a focused parameter domain $\mathcal{P}_{focus} \subset \mathcal{P}$ is determined, where a fine sampling of the parameter domain is used to obtain a highly accurate ROM solution.
This could induce a further speed up with respect to the full order model as the focused parameter domain may require less snapshot basis functions for an accurate ROM. 
A reference solution is computed for the real maximum of $\eta_{Lighthill}$ using the full order BEM.

\subsection{Stroke reconstruction}

To reconstruct the motion of the swimmer, we need to solve the complete time-dependent stroke.
Since the equation of motion are time-independent in the low Reynolds number regime,  
the time is simply indexing the geometrical variations $\boldsymbol{\mu}$ for the shape changes. 
Time thus behaves as any physical parameter for the model reduction purpose, which is in 
contrast to PDEs involving a time-derivative, where time is treated differently than other parameters. 
We assume complete knowledge of the geometry of the swimmer at all time instances and the solution of  training snapshots at the sample points.

\section{Numerical results}
\label{sec:results}
We present the numerical results obtained from coupling the boundary element method and reduced order modeling. 
Firstly, we present the convergence of the ROM to the full order solution with increasing dimension of the reduced model: the main focus is on the POD approximation \cite{Hesthaven2016}, 
which is tested on both applications in Section~\ref{sec:PODresults}.
In Section~\ref{sec:PODGreedy} we compare the POD and Greedy approach on the Eukaryotic-like swimmer, where the POD approximation turns out to be less accurate. 
Then we present two applications: the shape optimization of the bacterium-like swimmer, 
see Section~\ref{sec:2stepoptimization}, and the stroke reconstructions for the Eukaryotic-like swimmer, see Section~\ref{sec:strokereconstruction}.

\subsection{POD approximations of BEM solution}
\label{sec:PODresults}
In Sections~\ref{sec:split_rom} and~~\ref{sec:monolithic_rom} we introduced the split approach and the monolithic approach for projection-based ROMs for micro-swimmers. 
We compare the two approaches both for the bacterium-like and Eukaryotic-like swimmer, 
we consider both convergence to the full order solutions and online timings required to reach a prescribed accuracy. 
For this analysis we consider the maximum accuracy for the EIM and we let the number of modes gradually increase to see the convergence of the ROMs.
The split approach allows for 7 different projection spaces for the different tractions, we report only the analysis for the traction associated to the shape force $f_{shape}$ since 
it is the most demanding one in terms of accuracy and number of modes.

\subsubsection{Robotic bacterium}
\label{sec:PODbacterium}
\begin{figure}
\centering
\includegraphics[width=.4\textwidth]{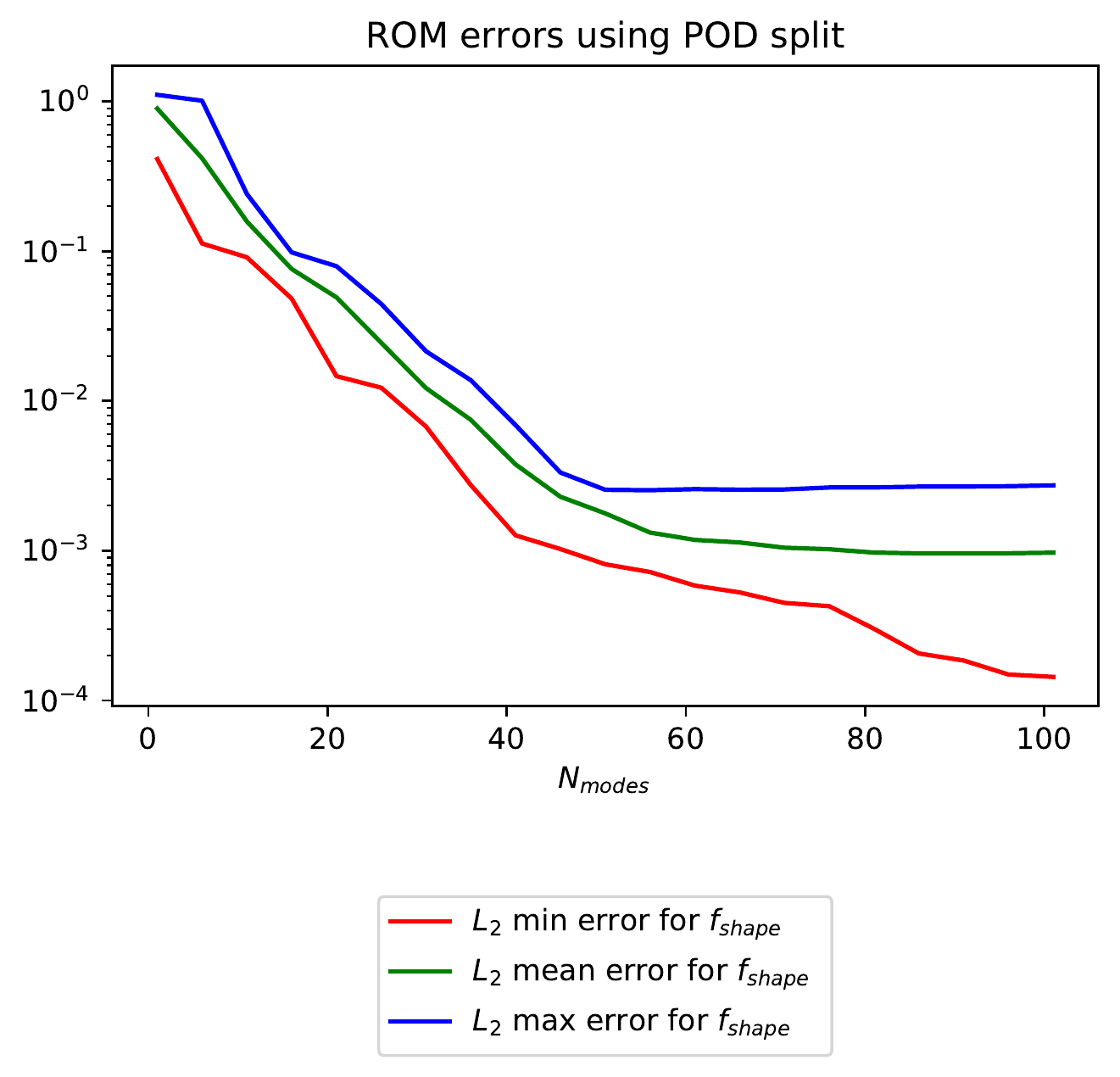} \ \includegraphics[width=.4\textwidth]{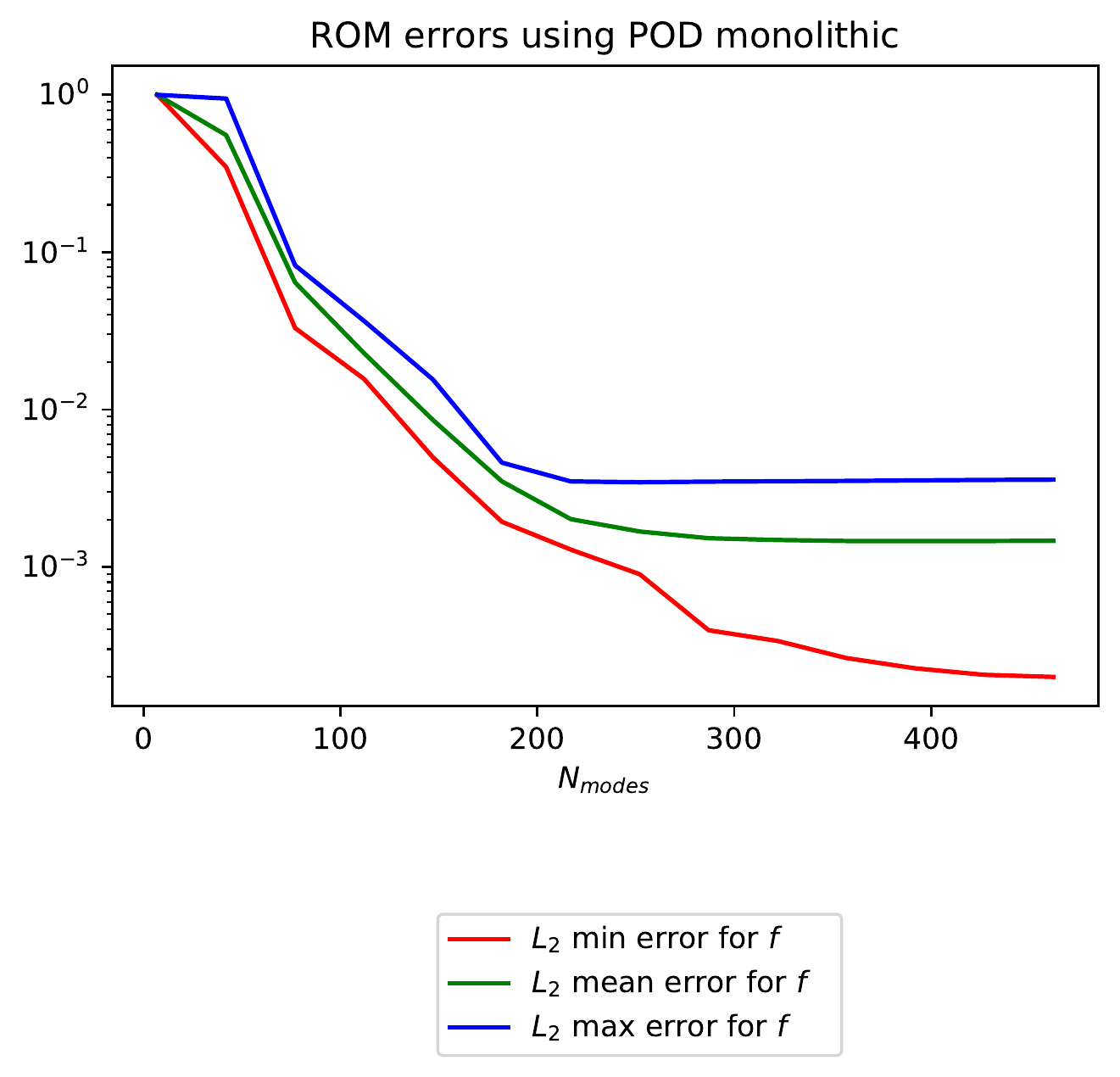}
\caption{ROM convergence analysis for the robotic-like swimmer computed for 10 random non training snapshots using POD. On the left the errors of the shape traction considering the split approach, on the right the final traction errors  using the monolithic approach. Red, green and blue represent the minimum, mean and maximum error respectively.}  
\label{PODerrorsBacterium}
\end{figure}

We compare split and monolithic POD approaches on the robotic-like bacterium presented in Section~\ref{subsec:micro_robot_models_bacterium}. 
Firstly we compare the convergence of the ROMs to the full order solution, in particular we analyze the errors relative to the traction vector $f$ in the two different models. 
Figure~\ref{PODerrorsBacterium} shows the convergence rate for the split approach on the left and for the monolithic approach on the right, red, green and blue depicts minimum, 
mean and maximum error computed over 10 randomly selected non training snapshots. 
We remark that the monolithic approach is built considering all the possible $N_{rigid} + 1$ tractions so this approach has an overall number of modes that is $N_{rigid} + 1$ times 
bigger than the split approach that requires $N_{rigid} + 1$ expansions instead. The convergence rates are very similar between the two different models, and they both reach a similar maximum 
accuracy between $10^{-5}$ and $10^{-3}$. If we increase further the number of modes a clear plateau emerges in the convergence and it is  mainly due to the errors introduced by the EIM. 
We tested the convergence of the ROM without considering the matrix approximations and we managed to achieve between four and five digits of 
accuracy without such an evident plateau.

In Table~\ref{table:summary_POD_bacterium} we compare the performances of the two POD models, we report the number of modes required to reach a prescribed accuracy of $5 \%$ and $0.5 \%$ both for the mean and maximum error. The split approach requires less number of modes than the monolithic one, however we remind that the latter approach requires a single resolution to get the final solution while the split approach solves $7$ different linear systems. From this analysis we see that the split approach has a slight advantage with respect to the monolithic one. We believe that this is due to a better memory handling of the split approach, in fact a deeper analysis reveals that the monolithic approach requires more time to compute the BEM entries to be used in the EIM, see Section~\ref{sec:DEIM}. The two approaches require the same entries to be computed but the reduced number of modes of the split approach makes available more memory to the BEM computation and this results in this slight speed up. Apart this minor difference we believe the two approaches to be equivalent.

\begin{table}
\centering
\begin{tabular}{| l | l | l | l | l |}  
\hline
\hline
  & Split mean & Monol mean   & Split max & Monol max \\
\hline
err. less $5 \%$  & 21 & 77  & 26 & 112\\
err. less $0.5 \%$  & 41 & 182   & 46 & 182\\
online timing  & 1.322 & 1.446   & 1.371 & 1.446\\
\hline
\hline
\end{tabular}
\caption{Summary of the ROM performance indicators for the POD approximations on the bacterium-like test case. Shown are the smallest basis sizes, where a mean approximation error below $5 \%$ and $0.5 \%$ is attained. The online timing is referred to the approximation needed to have $0.5 \%$ accuracy}
\label{table:summary_POD_bacterium}
\end{table} 

\subsubsection{Eukaryotic swimmer}
\label{sec:PODchlamy}
\begin{figure}
\centering
\includegraphics[width=.4\textwidth]{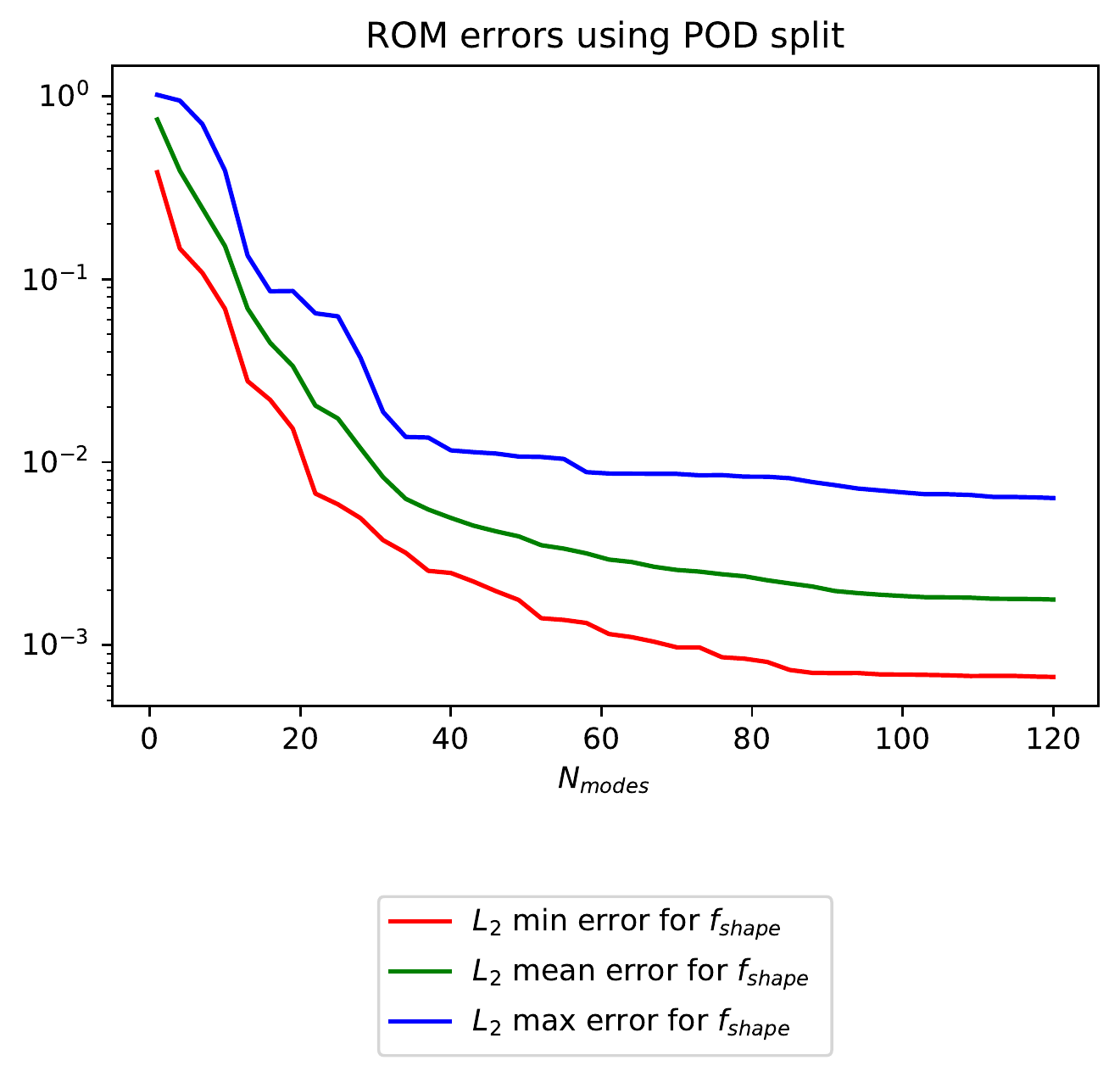} \ \includegraphics[width=.4\textwidth]{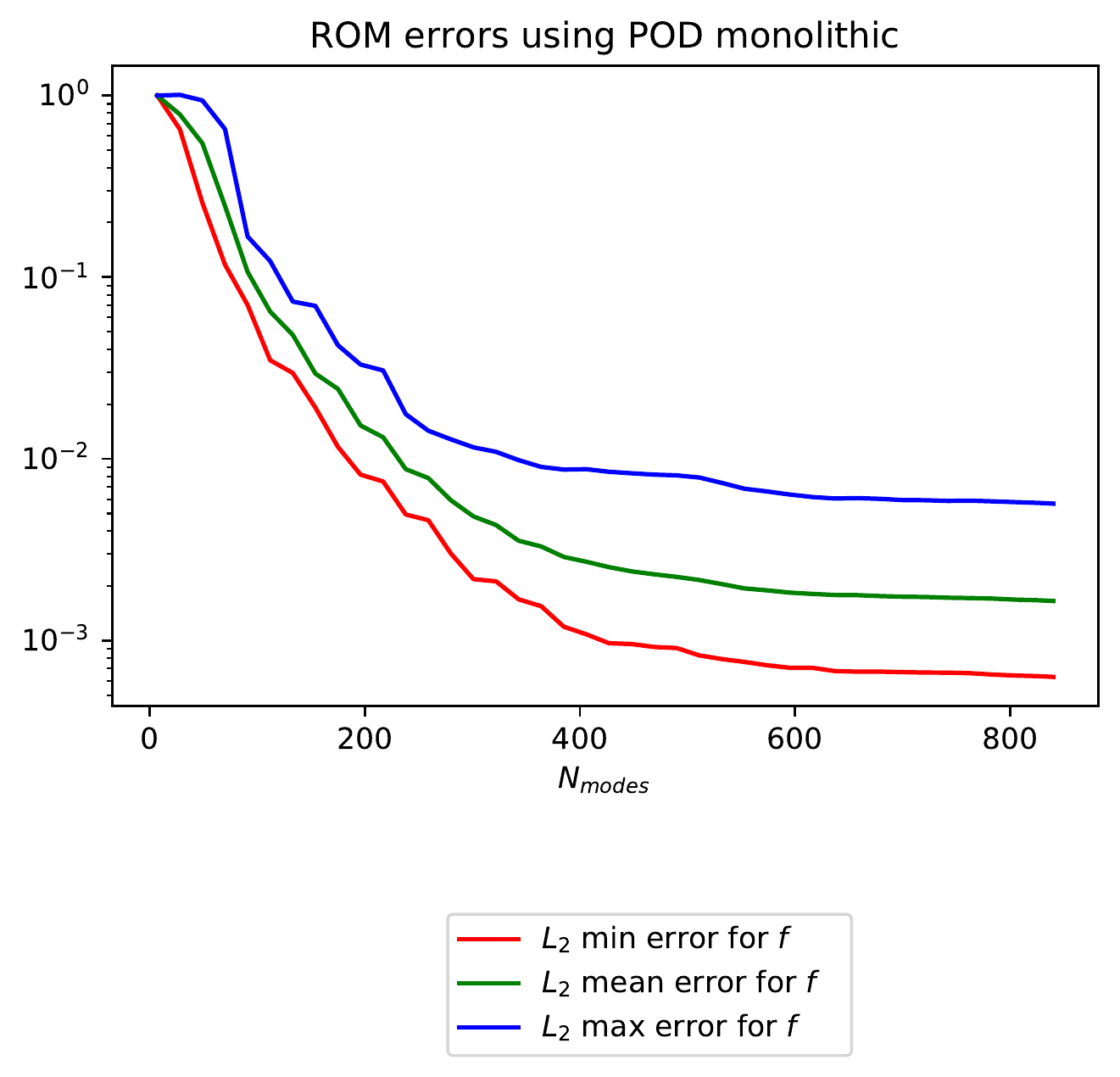}
\caption{ROM convergence analysis for the Eukaryotic-like swimmer computed for 10 random non training snapshots using POD. On the left the errors of the shape traction considering the split approach, on the right the final traction errors using the monolithic approach. Red, green and blue represent the minimum, mean and maximum error respectively.}  \label{PODerrorsSplitChlamydomonas}
\end{figure}

We compare the convergence to the full order model and the performances of the two POD approaches on the Eukaryotic-like swimmer presented in Section~\ref{subsec:micro_robot_models_chlamydomonas}. Figure~\ref{PODerrorsSplitChlamydomonas} compares the convergence of the split approach (on the left) and the monolithic approach (on the right). Red, green and red represent minimum, mean and maximum errors computed on 10 randomly chosen non training snapshots. 
We see that the two ROMs present the same convergence to the full order solution, the error in both cases converges between $10^{-3}$ and $10^{-2}$. We don't see any clear plateau in this test-case and this is due to the increase complexity in the geometry for the Eukaryotic like movement we consider. For this reason the EIM plateau of Section~\ref{sec:PODbacterium} is not reacehd.

Table~\ref{table:summary_POD_chlamy} compares the performances of the two approaches considering the number of modes required to reach $5\%$ and $0.7\%$ accuracy. As expected, the split approach requires less modes to reach the target accuracy but we note again that the major part of the online time is required to reconstruct the matrices using EIM, and specifically it is needed to compute the BEM entries. If the number of modes remains low the timings are the same while for the last test case when we require $553$ modes we notice a difference in the timings. We believe that this phenomenon confirms our speculation about the worse memory handling of the code when the number of modes increases leaving less memory available to the BEM computations. Apart from this non-linear minor effect we see that the two methods are almost equivalent both in terms of convergence and performances. 
\begin{table}
\centering
\begin{tabular}{| l | l | l | l | l |}  
\hline
\hline
  & Split mean & Monol mean  & Split max & Monol max \\
\hline
err. less $5 \%$  & 16  & 133  & 28 & 175\\
err. less $0.7 \%$  & 34  & 280  & 100 & 553\\
online timing  & 0.6637  & 0.6631  & 0.6792 & 1.107\\
\hline
\hline
\end{tabular}
\caption{Summary of the ROM performance indicators for the POD approximations for the Eukaryotic-like swimmers. Shown are the smallest basis sizes, where a mean approximation error below $5 \%$ and $0.7 \%$ is attained. The online timing is referred to the approximation needed to have $0.7 \%$ accuracy}
\label{table:summary_POD_chlamy}
\end{table} 

\subsection{Greedy and POD: comparison of approximation accuracy}
\label{sec:PODGreedy}
\begin{figure}
\centering
\includegraphics[width=.4\textwidth]{chlamy_rom_errors_random_10.pdf} \ \includegraphics[width=.4\textwidth]{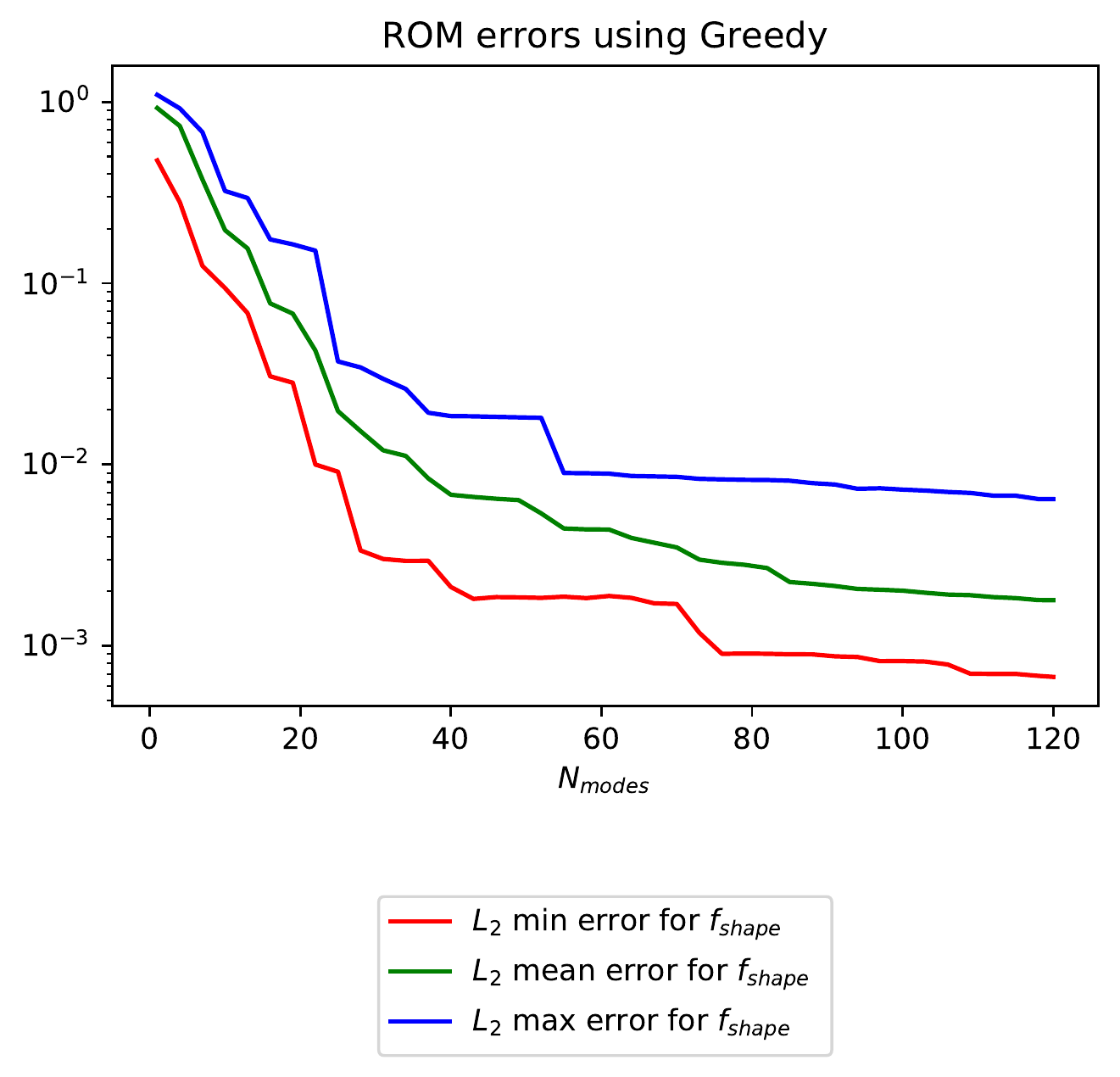}
\caption{Comparison of ROM convergence  for the Eukaryotic-like swimmer computed for 10 random non training snapshots using POD or Greedy approach. On the left the POD errors of the shape traction, on the right the shape traction errors using Greedy. Red, green and blue represent the minimum, mean and maximum error respectively.}  
\label{GreedyerrorsSplitChlamydomonas}
\end{figure}

\begin{figure}
\centering
\includegraphics[width=.4\textwidth]{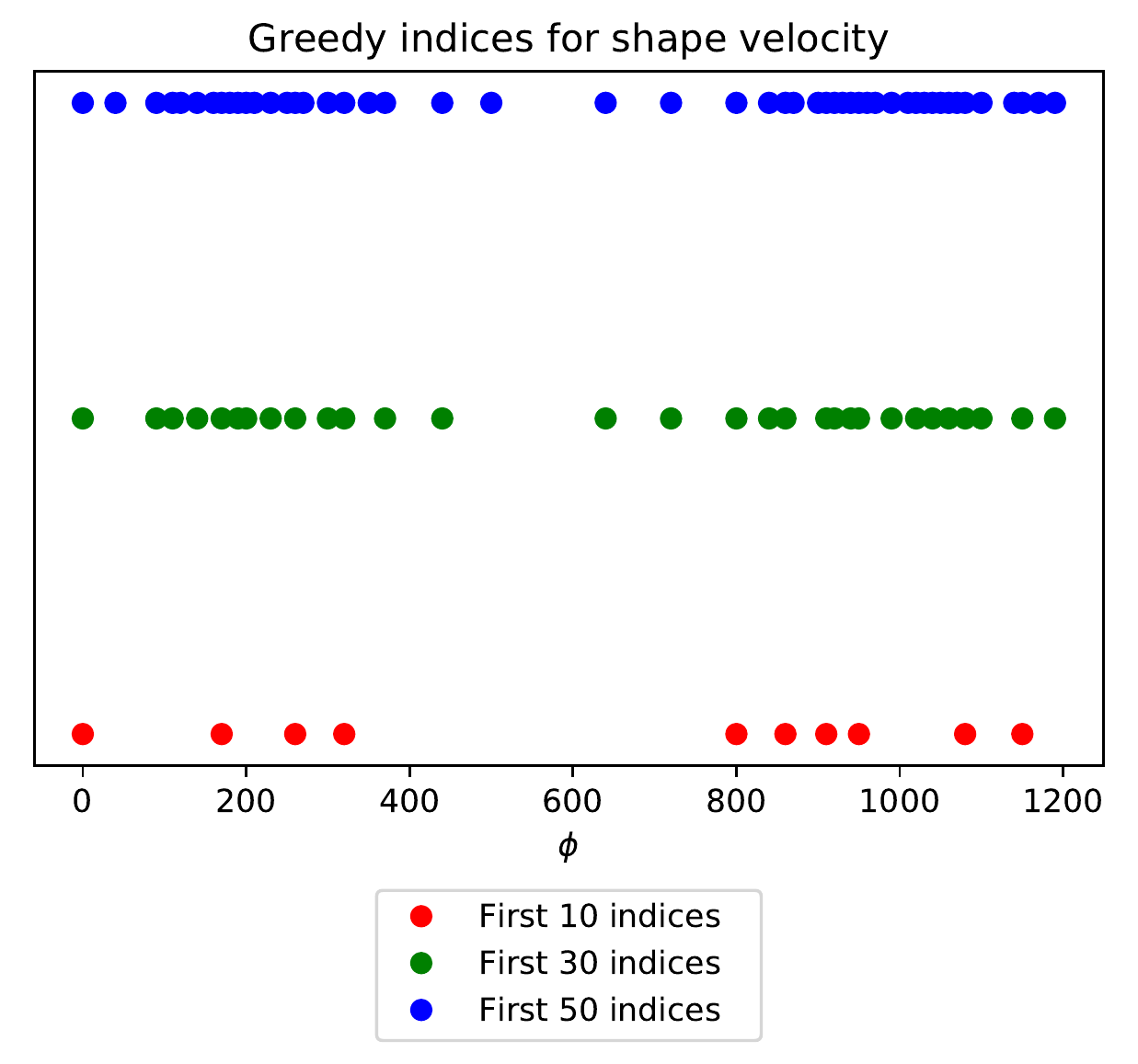} \ \includegraphics[width=.4\textwidth]{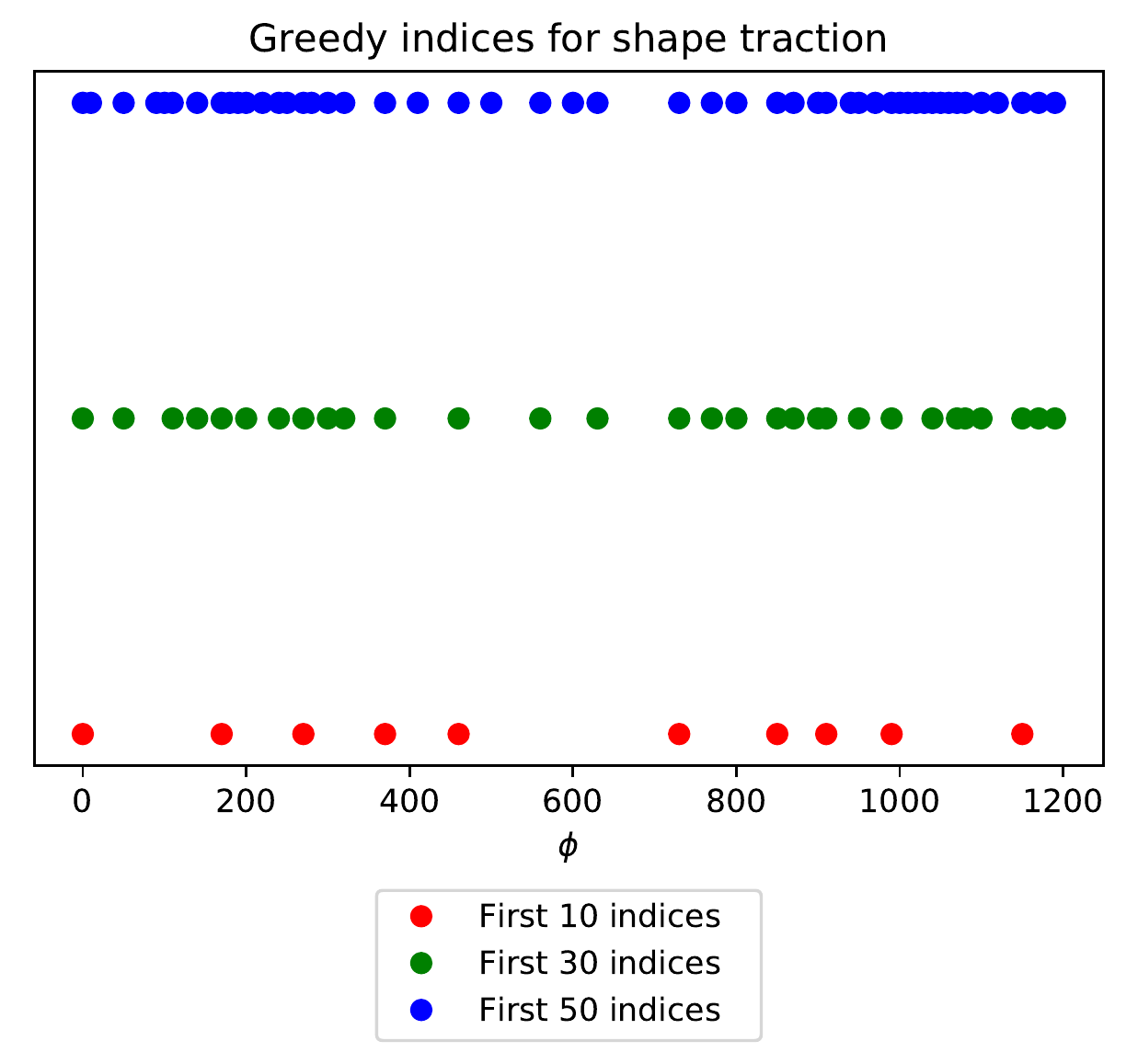}
\caption{Greedy snapshot selection for shape velocity (on the left) and shape traction (on the right) approximations. Red, Greed and blue shows the first $10, 30, 50$ selected modes.}  
\label{Greedyindices}
\end{figure}

\begin{figure}
\centering
\includegraphics[width=.45\textwidth]{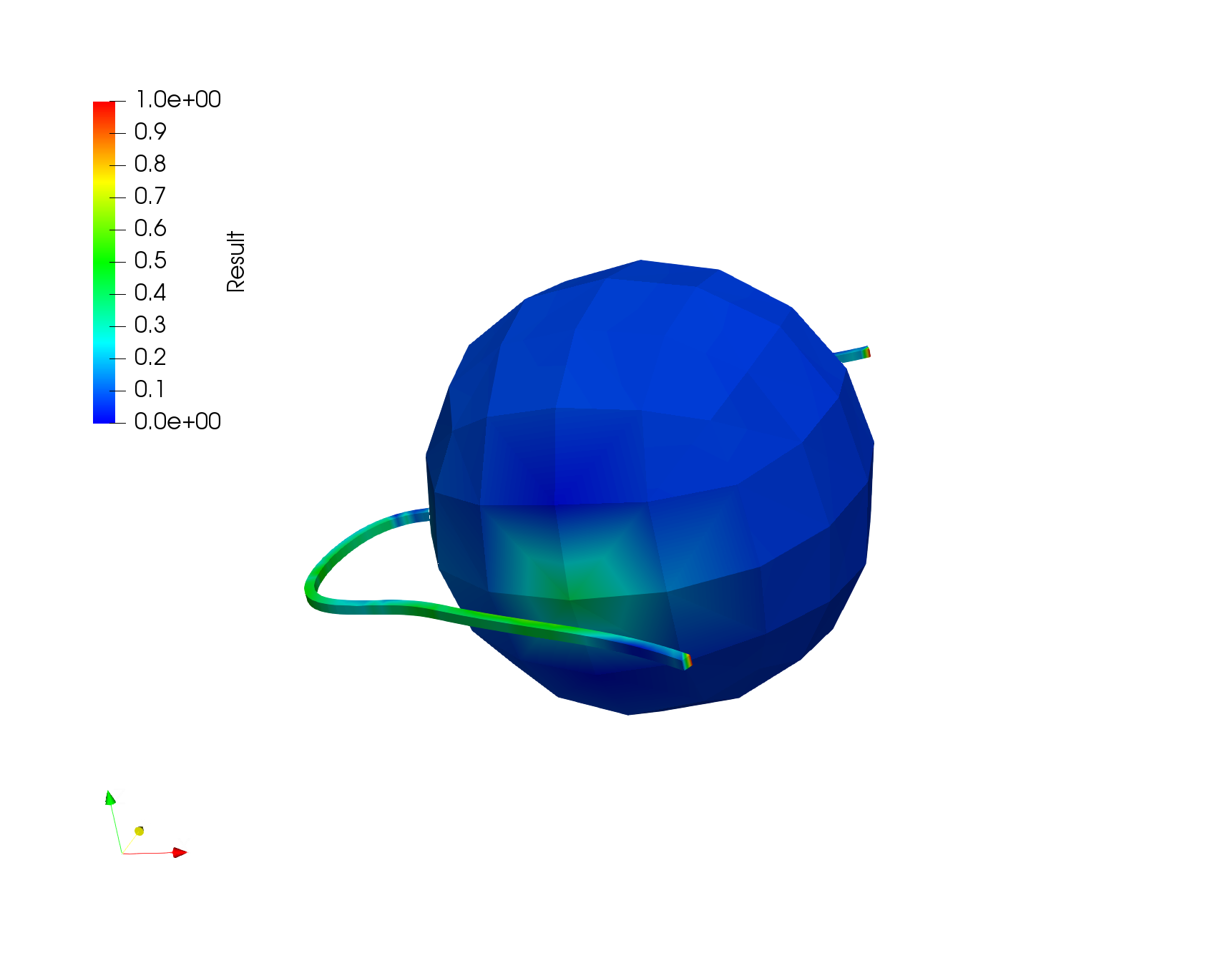} \includegraphics[width=.45\textwidth]{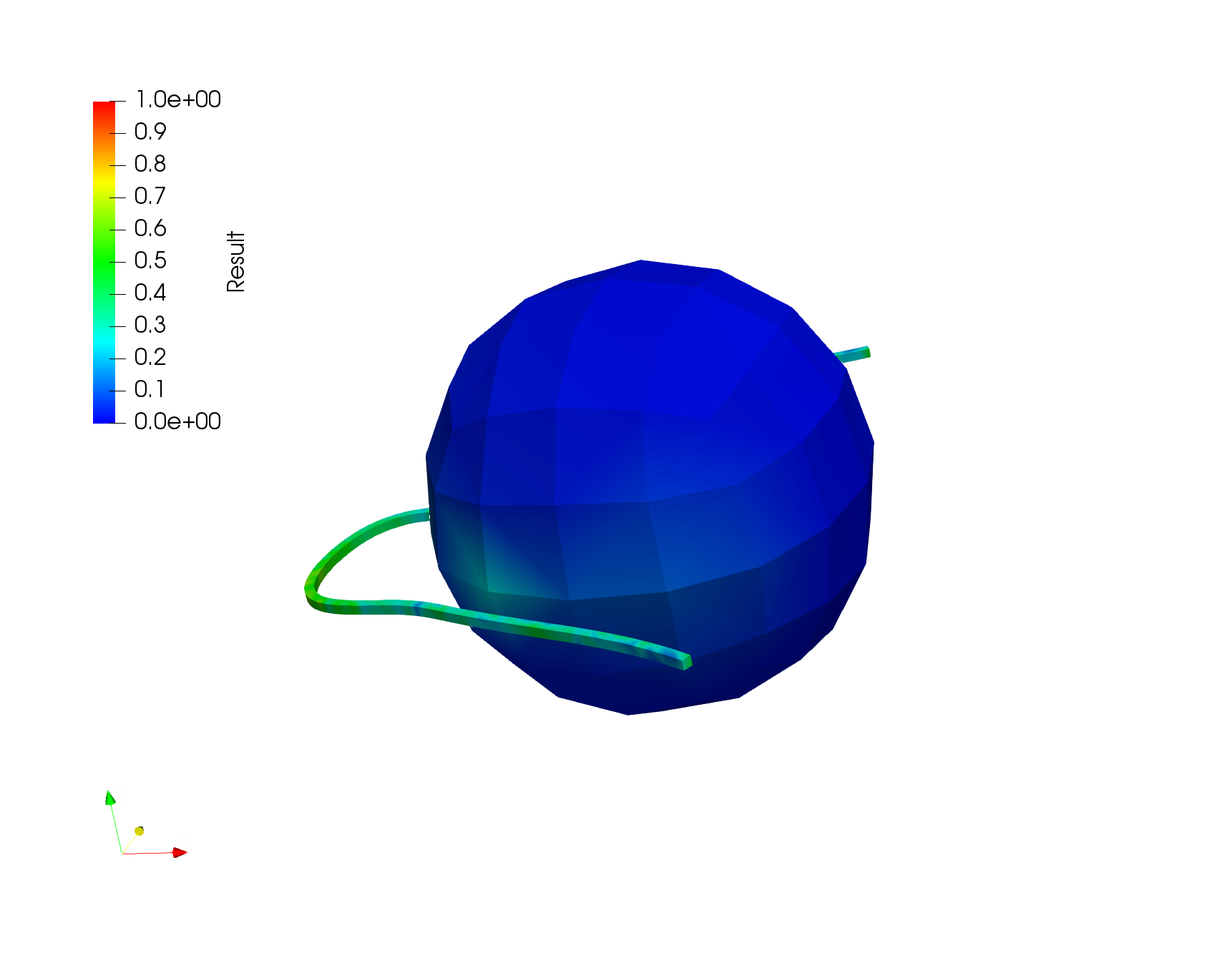}
\caption{Comparison of the traction error magnitude for the Eukaryotic swimmer when the flagella are very near the cell body. On the left the reconstructed solution using $25 \%$ of the available modes, on the right the most accurate reduced order solution. The scale is the same on the two plots and it is set so that 1 represent the maximum error of the two approximations.}  
\label{ReconstructedTractionMagnitude}
\end{figure}

In Section~\ref{sec:PODresults} we saw that the two POD approaches are mostly equivalent on both the test-cases presented in Section~\ref{sec:micro_robot_models}. However we noticed that the geometry shape change of the Eukaryotic-like swimmer, represented in Figure~\ref{chlamy}, is more complicated than the rotation of the bacterium-like one.  This induces a worse convergence to the full order solution. For this reason we compare the split POD approach to the Greedy approach introduced in Section~\ref{sec:greedy_split} on this test-case.
Figure~\ref{GreedyerrorsSplitChlamydomonas} compares the convergence of the POD approach (on the left) and the Greedy approach (on the right), red, green and blue represent minimum, mean and maximum errors computed on 10 randomly chosen non training snapshots. 
We see that the two approaches have a very similar convergence to the exact solution, we notice that the POD approach is smoother but they both reach the same final accuracy. 
Table~\ref{table:summary_Greedy} compares the performances of the two approaches and we see that the POD and the Greedy are equivalent even from this point of view. 
\begin{table}
\centering
\begin{tabular}{| l | l | l | l | l |}  
\hline
\hline
  & POD mean & Greedy  mean & POD max & Greedy max\\
\hline
err. less $5 \%$  & 16 & 22 & 28 & 25\\
err. less $0.7 \%$  & 34 & 40 & 100 & 109 \\
online timing  & 0.6637 & 0.6547 & 0.6792 & 0.6814\\
\hline
\hline
\end{tabular}
\caption{Summary of the ROM performance indicators for the POD and Greedy approximation on the Eukaryotic-like test-case. Shown are the smallest basis sizes, where a mean approximation error below $5 \%$ and $0.7 \%$ is attained. The online timing is referred to the approximation needed to have $0.7 \%$ accuracy}
\label{table:summary_Greedy}
\end{table} 
To better understand the Greedy selection procedure we represent in Figure~\ref{Greedyindices} the first indices selected to approximate the shape velocity $v$ (on the left) and the corresponding traction force $f_{shape}$ (on the right). 
The indices cluster around two different points: the first point is located at $\phi \sim 200$ which corresponds to the transition between backward and forward motion, the second one is at $\phi \sim 900$ which corresponds to the two flagella moving very close to the body. Since inertia is negligible the transition between backward and forward motion is immediate as soon as the flagellar beat allows it and the Greedy procedure selects frame near the transition. When the flagella are close to the body there is a lot of interaction between different body parts and the Greedy procedure selects frames in this region to have a better approximation. To better understand this scenario we compare the magnitude of the traction when flagella are very near the body in Figure~\ref{ReconstructedTractionMagnitude}. 
Even if the differences are  moderate the solution with one quarter of the possible Greedy modes depicted on the left is not able to recover the traction pattern of the reference solution on the right. Consequently, the Greedy algorithm tries to increment the accuracy by selecting more modes for this situation, the plot at the centre shows the most accurate possible solution which has a slight increase in the traction representation.

From our analysis on the Eukaryotic-like swimmer we conclude that Greedy and POD have the same convergence and timings. The results of Sections~\ref{sec:PODbacterium} and~\ref{sec:PODchlamy} show the same behavior for all the POD on both application so we expect Greedy and POD to behave similarly with the monolithic approach and the bacterium-like application. This grants flexibility for further use cases as
either POD / Greedy and split / monolithic approaches can achieve speed-ups of several magnitudes.

\subsection{Two step shape optimization of robotic bacterium}
\label{sec:2stepoptimization}
In this Section we study the actual shape optimization of the robotic micro-swimmer introduced in Section~\ref{subsec:micro_robot_models_bacterium}. As performance measure we consider the energetic efficiency $\eta_{Lighthill}$ introduced in~\eqref{eta_lighthill}. We find the optimal value of the parameter $\boldsymbol{\mu} = (N_\lambda, R_{head})$ in the parameter domain $\mathcal{P} = [0.4, 4.0]^2$. We only use one of the possible ROMs we presented since they proved to have very similar accuracy and performances. 
A common way to study and optimize the performance of the considered swimmer is to study the hydrodynamics of the separate component (head and flagellum) to infer the swimming behavior of the complete swimmer~\cite{Purcell1977,Purcell1997}. For a complete analysis of this ``additive approach'' (AA) the reader is referred to~\cite{Giuliani2018SoRo}. 

\begin{figure}
\centering
\includegraphics[width=.4\textwidth]{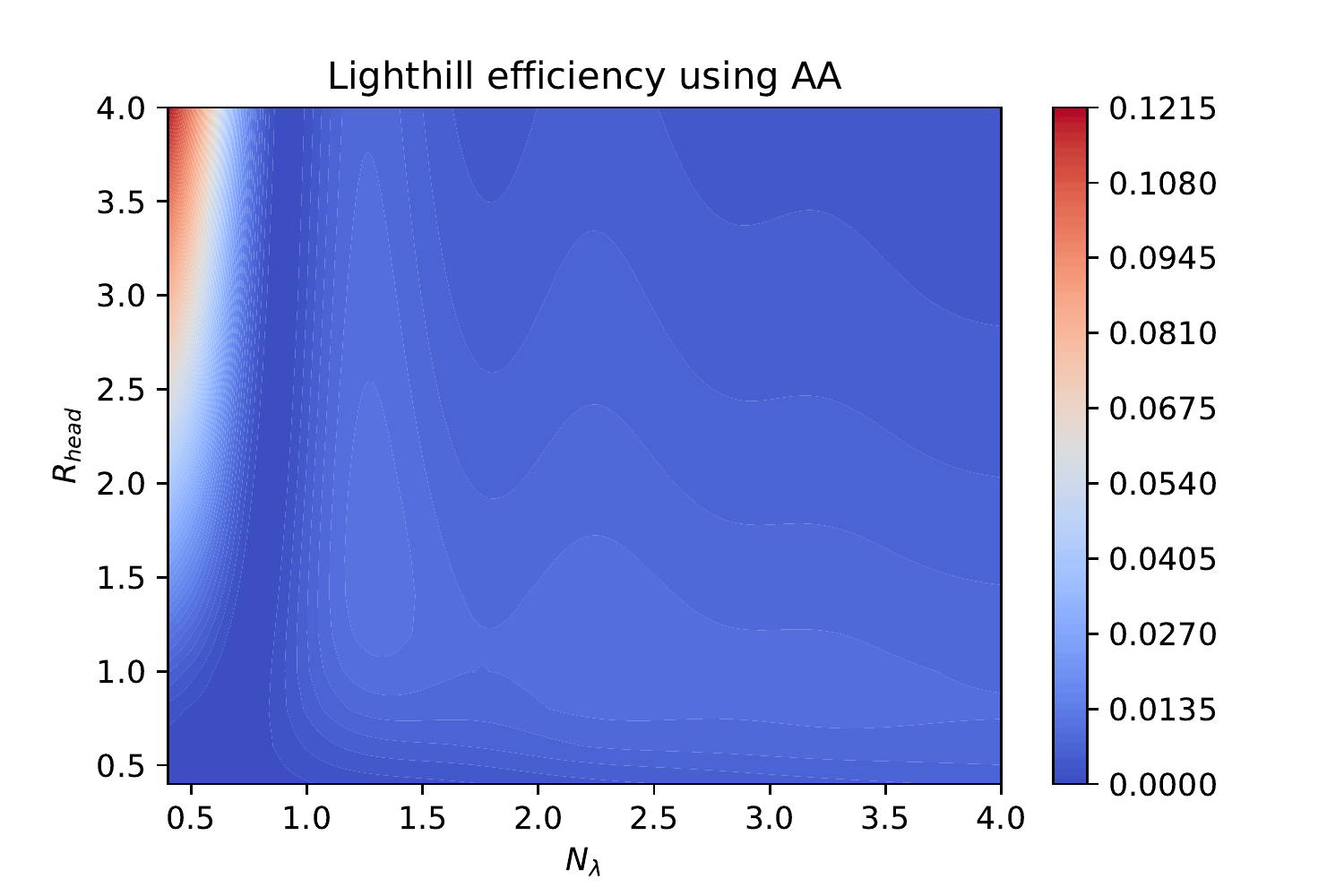} \includegraphics[width=.4\textwidth]{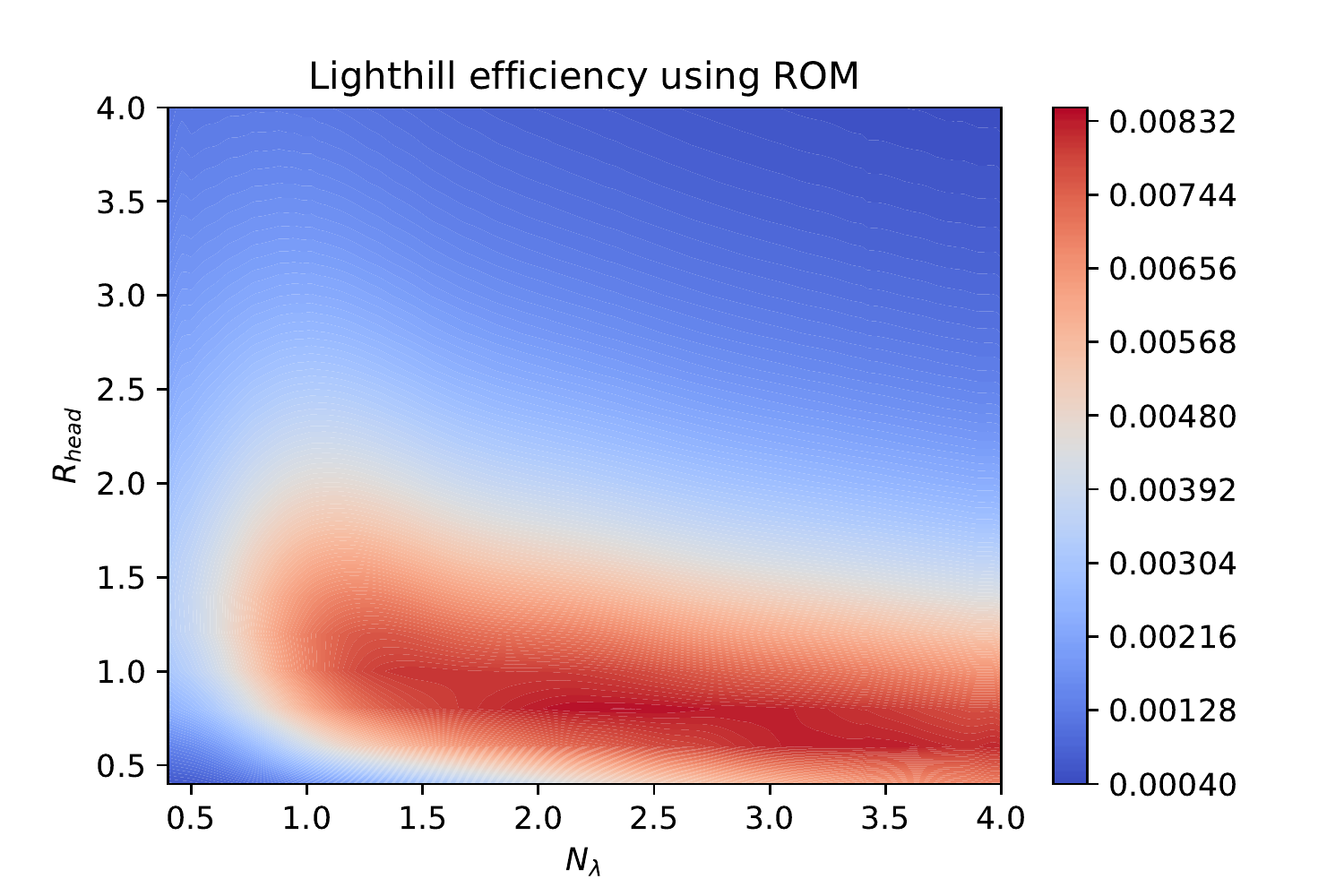}
\caption{Lighthill efficiency~\eqref{eta_lighthill} analysis $N_{\lambda} = 0.4 - 4.0$, $R_{head} = 0.4 - 4.0$. AA on the left, ROM on the right.}
\label{LighthillEfficiencyCoarse}
\end{figure}

In Figure~\ref{LighthillEfficiencyCoarse} we compare $\eta_{Lighthill}$ analysis using AA (on the left) and ROM (on the right) letting the radius of the head $R_{head}$ to vary with a step of $0.2$ and the number of turns for the tail with a step $0.02$.
AA shows the optimum ($\eta_{Lighthill}=0.12138$) for $N_{\lambda}=0.4, R=4.0$. Considering the analysis of~\cite{Giuliani2018SoRo} we believe this optimal value to be erroneous since a big head radius creates a considerable wake effect, which is completely neglected by AA, worsening the performances of the real complete swimmer. 
ROM shows instead an optimum ($\eta_{Lighthill}=0.008406$) for $N_{\lambda}=2.36, R=0.8$ and we see that also AA shows a local maximum close to $N_\lambda = 2.4$ and $R = 0.7$. 
We compute the real maximum for $\eta_{Lighthill}$ using the full order BEM on the same parameter space and we obtain the optimum at $N_{\lambda}=2.38, R=0.8$ with a relative error for the efficiency of $1.41 \%$ for the ROM.  

\begin{figure}
\centering
\includegraphics[width=.4\textwidth]{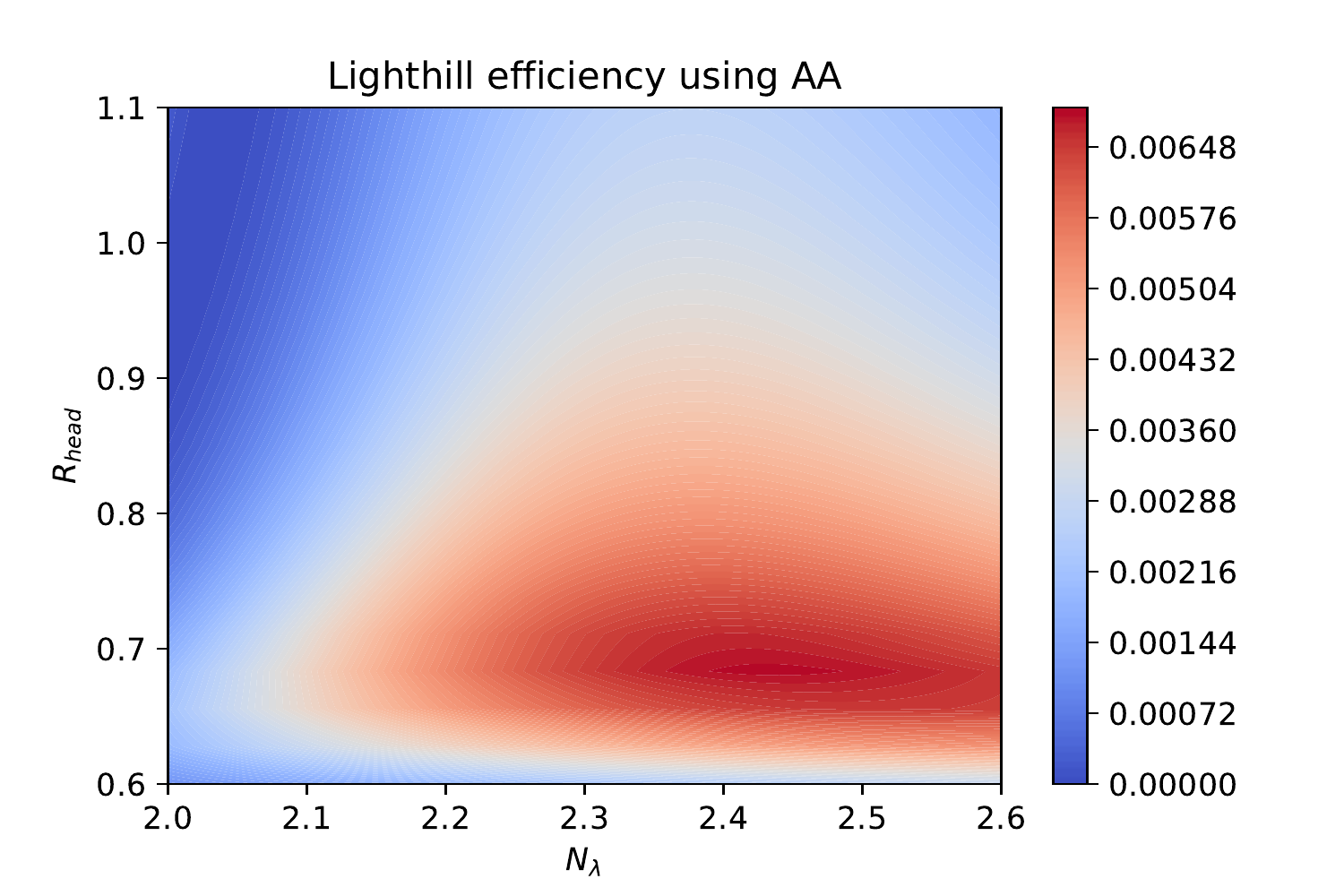} \includegraphics[width=.4\textwidth]{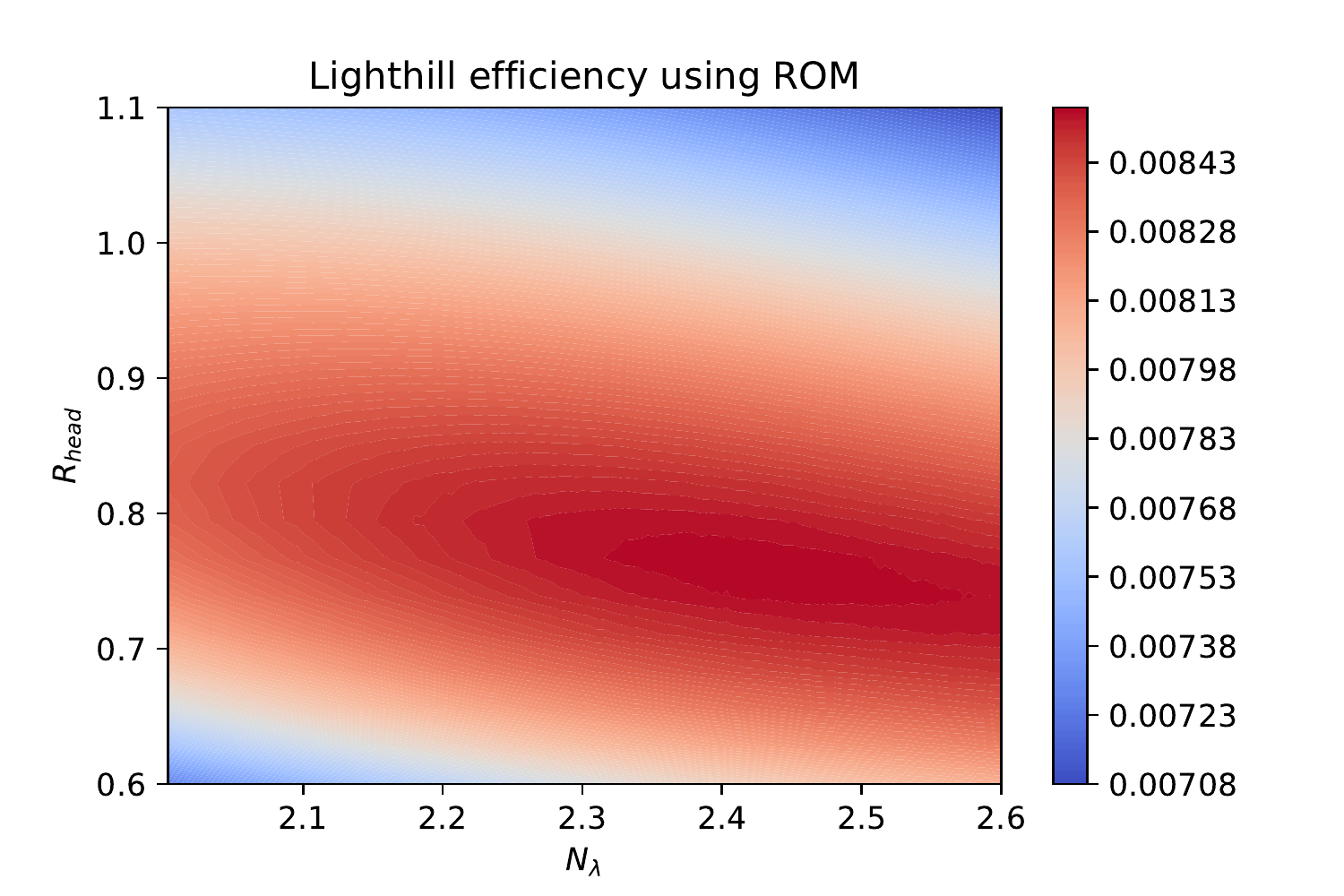}
\caption{Lighthill efficiency~\eqref{eta_lighthill} analysis  $N_{\lambda} = 2.0 - 2.6$, $R_{head} = 0.6 - 1.1$. AA on the left, ROM on the right.}
\label{LighthillEfficiencyFine}
\end{figure}

To increase the accuracy of our predictions we apply a second optimization step using a ROM build considering a parameter grid centered on the maximum value region $\mathcal{P}_{focus} \subset \mathcal{P}$, namely we consider a new set of training $\boldsymbol{\mu} = (N_\lambda, R_{head})$ in the parameter domain $\mathcal{P}_{focus} = [0.6, 1.1] \times [2.0, 2.6]$. Figure~\ref{LighthillEfficiencyFine} shows the results on the fine grid both for AA (on the left) and ROM (on the right).
AA predicts an optimum ($\eta_{Lighthill}=0.00684121$) for  $N_{\lambda}=2.4333, R_{head}=0.6833$, ROM shows the optimum ($\eta_{Lighthill}=0.0085464$) for  $N_{\lambda}=2.42, R_{head}=0.7667$, while the full order BEM displays the optimum $\eta_{Lighthill}=0.00854549$) for  $N_{\lambda}=2.4067, R_{head}=0.7667$. 
The simplified AA is not able to recover a good result, especially for what concerns the head radius, while a proper ROM can be effectively used to compute a configuration which is very close to the optimal one for a robotic micro-swimmer.
Table~\ref{error_sum} sums up the various predictions for the maximum $\eta_{Lighthill}$ for all the different scenarios and we clearly see that the ROM is able to recover very good solutions for the optimal shape.
\begin{table}
\centering
\begin{tabular}{| l | l | l | l | l |}  
\hline
\hline
 Method & $R_{head}$ & $N_{turns}$ & $\eta_{Lighthill}$ & Error \\
\hline
BEM coarse  & $0.8$ & $2.38$ &  $0.008526$  &  \\
AA coarse & $4.0$ & $0.4$ &  $0.12138$ &  $1323 \%$\\
ROM coarse & $0.8$ & $2.36$ &  $0.008406$ &  $1.41\%$ \\
BEM fine  & $0.7667$ & $2.4067$ &  $0.00854549$ &   \\ 
AA fine & $0.6833$ & $2.4333$ &  $0.00684121$ &  $19.94 \%$ \\
ROM fine & $0.7667$ & $2.42$ &  $0.0085464$ &  $0.0104 \%$ \\
\hline
\hline
\end{tabular}
\caption{$\eta_{Lighthill}$ predictions and errors using both additive approach and ROM.}
\label{error_sum}
\end{table}

\subsection{Stroke reconstruction of a micro-swimmer}
\label{sec:strokereconstruction}
Another possible application of the ROM procedures described in this work is the reconstruction of the rigid body velocities of the swimmer during the complete stroke starting from some precomputed snapshots. We apply this procedure to the Eukaryotic-like swimmer presented in Section~\ref{subsec:micro_robot_models_chlamydomonas}. We have already shown that all the presented ROMs leads to very similar results, so in this Section we only use one of them, in particular we use the POD split approach.

\begin{figure}
\centering
\includegraphics[align=t,width=.4\textwidth]{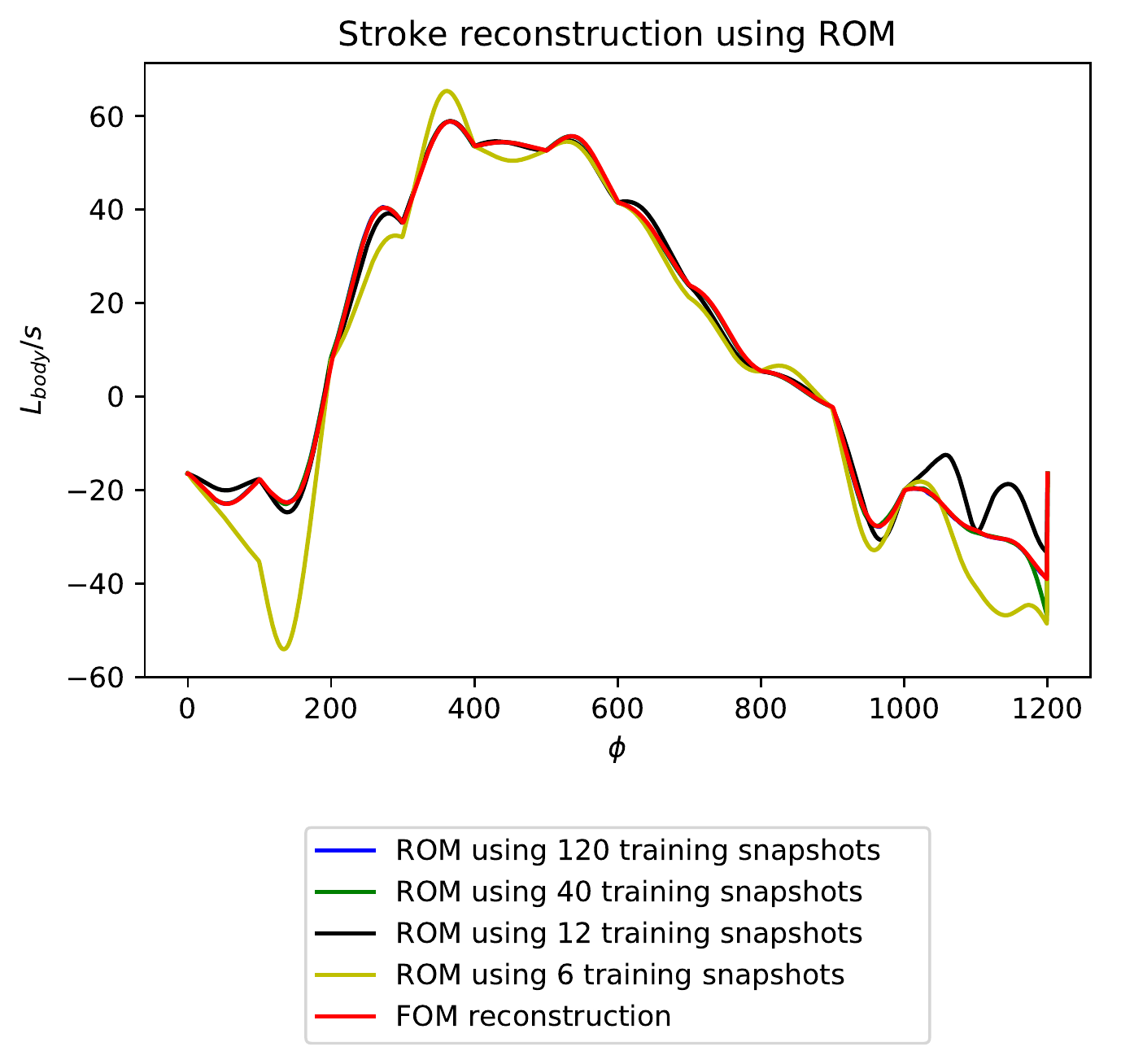} \ \includegraphics[align=t,width=.4\textwidth]{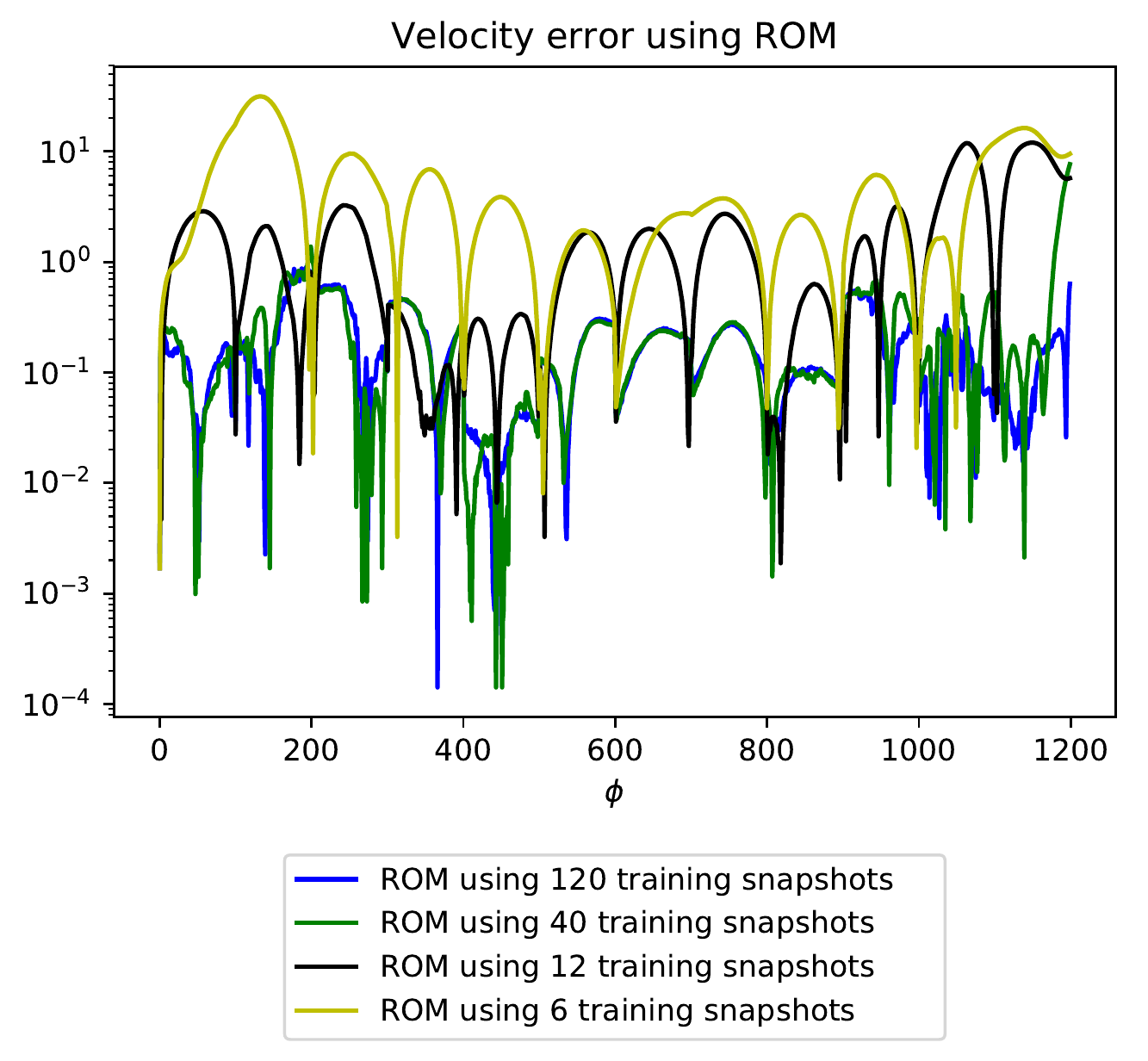}
\caption{Stroke reconstruction of an Eukaryotic-like swimmer. On the left we plot the velocity along the longitudinal axis. On the right we report the error of the ROM approximations with respect to the full order model. Yellow, black, green and blue lines define ROM approximations built considering an increasing number of training snapshots ($6,12,40,120$), the red line represents the reference solution obtained with the full order BEM.}  
\label{stroke_reconstruction}
\end{figure}

On the left of Figure~\ref{stroke_reconstruction} we use ROM to reconstruct the rigid velocities during the complete stroke starting from the knowledge of the system at the training snapshots. The artificial stroke that we want to reconstruct consists of $1200$ different frames, we consider an increasing number of equi-spaced training snapshots ($N_{training} = 6, 12, 40, 120$) and we see that the solution converges to the full order solution increasing the complexity of the ROM. 
 ROM has proved to be an effective tool to reconstruct the complete stroke starting from a limited number of training snapshots. The right plot of Figure~\ref{stroke_reconstruction} clearly shows the convergence of the stroke to the proper one as we increase the number of training snapshots.

\section{Conclusions}
\label{sec:conclusions}

This work shows that ROM enhanced BEM simulations of microswimmers can achieve accurate results while significantly reducing the runtime.
The many-query setting, which is typical to ROM, also holds when resolving models depending only on time as a single parameter.
This is because inertia, i.e., a time-derivative, is not present and the fluid ``reaction'' is quasi instantaneous.
This widens the applicability of ROM methods to microswimmers significantly, as each use-case can be treated separately in a ROM or
several use-cases can be further parametrized so that ROM can be applied to a whole family of microswimmers.
Another feature of the presented approach is the geometric flexibility, which applies EIM to online-generated mesh data. 
Using the high-order mesh generation in the online phase allows significantly more geometric flexibility than the restriction to affine geometry transformations, and it allows to capture also large and non-linear deformations.

Both tested formalisms of the BEM systems (i.e., split and monolithic approach) and standard reduced basis ROM algorithms (i.e., POD and greedy) 
are capable to mitigate the computational costs with an acceptable approximation accuracy. 
Each approach reaches engineering standards of on average three digits of accuracy.
The methodology is applied on two different test-cases, and the maximum reduced order solution error is below $1\%$ in both models, the robotic and eukaryotic.

To be able to further increase the accuracy of the presented methodology a meshing software which
provides more than single digit accuracy is necessary.
With the presented approach of EIM, projection, solving the reduced order system and approximating 
the output quantity, on average 3 to 4 digits of accuracy on the output quantity can be achieved.

\section*{Acknowledgments}
\noindent
This work has been supported by the FP7 ERC Advanced Grant 2013 project 340685 ``MicroMotility" and H2020 ERC Consolidator Grant 2015 AROMA-CFD
project 681447 ``Advanced Reduced Order Methods with Applications in
Computational Fluid Dynamics". Nicola Giuliani thanks Dr. Giovanni Stabile, for useful discussions.
The support of INDAM GNCS is also acknowledged.

\bibliographystyle{plain}
\bibliography{bibliography}

\end{document}